\documentclass[12pt]{article}
\usepackage{amsmath, graphicx, amsfonts,amssymb, calrsfs}
 
\usepackage{color}

\def\bbR{\mathbb{R}}
\def\r{\rho}
\def\s{\sigma}
\def\a{\alpha}

\def\cK{\mathcal{K}}
\def\cV{\mathcal{V}}

\def\L{\Lambda}

\def\cF{\mathcal{F}}
\def\o{\omega}
\def\ball{B^n_2}
\def\polar{K^\circ}
\def\l{\lambda}
\def\cS{\mathcal{S}}

\def\cL{\mathbf{L}}

\def\bK{\mathbf{K}}
\def\be{\begin{equation}}
\def\ee{\end{equation}}
\def\bea{\begin{eqnarray}}
\def\eea{\end{eqnarray}}
\def\bt{\begin{theorem}}
\def\et{\end{theorem}}
\def\bl{\begin{lemma}}
\def\bQ{\mathbf{Q}}
\def\el{\end{lemma}}
\def\br{\begin{remark}}
\def\er{\end{remark}}
\def\bc{\begin{corollary}}
\def\ec{\end{corollary}}
\def\bd{\begin{definition}}
\def\ed{\end{definition}}
\def\bp{\begin{proposition}}
\def\ep{\end{proposition}}
\newtheorem{theorem}{Theorem}[section]
\newtheorem{lemma}{Lemma}[section]
\newtheorem{remark}{Remark}[section]
\newtheorem{proposition}{Proposition}[section]
\newtheorem{corollary}{Corollary}[section]

\newtheorem{definition}{Definition}[section]
\begin{document}
\title{The mixed $L_p$ geominimal surface areas for multiple convex  bodies
}

\author{Deping Ye, Baocheng Zhu  and Jiazu Zhou}
\date{}
\maketitle

\begin{abstract}
In this paper, we introduce several mixed $L_p$ geominimal surface areas for multiple convex bodies for all $-n\neq p\in \bbR$. Our definitions are motivated from an equivalent formula  for the mixed $p$-affine surface area. Some properties, such as the affine invariance,  for these mixed $L_p$ geominimal surface areas are proved. Related inequalities, such as, Alexander-Fenchel type inequality, Santal\'{o} style inequality, affine isoperimetric inequalities, and cyclic inequalities are established. Moreover, we also study some properties and inequalities for the $i$-th mixed $L_p$ geominimal surface areas for two convex bodies.  \vskip 2mm
2010 Mathematics Subject Classification: 52A20, 53A15
\end{abstract}

\section{Introduction}
The combination of the Minkowski sum and the volume naturally leads to the mixed volume for multiple convex bodies (i.e., convex compact subsets in $\bbR^n$ with nonempty interior), which is now the very core of the Brunn-Minkowski theory of convex bodies. Numerous widely-studied functionals on convex bodies, e.g., the volume and the surface area, are special cases of the mixed volume. The mixed volume has many nice properties and important inequalities which are fundamental in applications. For instance, the Alexander-Fenchel inequality related to the mixed volume is one of the most important inequalities in convex geometry. Many fundamental geometric inequalities, such as, the Minkowski's first inequality and the Brunn-Minkowski inequality for convex bodies, can be derived from the Alexander-Fenchel inequality. Readers are referred to \cite{Schn} for more details and more references regarding the mixed volume and related properties. 

The rapidly developing $L_p$ Brunn-Minkowski theory of convex bodies is a natural extension of the Brunn-Minkowski theory. Such an extension arises from the combination of the volume and the Firey $p$-sum of convex bodies for $p\geq 1$ introduced by Firey in \cite{Firey} about 50 years ago.  It has been thought that the  $L_p$ Brunn-Minkowski theory of convex bodies has the $L_p$ affine surface area as its core. The $L_p$ affine surface area was introduced by Lutwak in his seminal paper \cite{Lu96} about 70 years after the classical affine surface area (i.e., for $p=1$) was first introduced by Blaschke \cite{Blas} in 1923. One of the most remarkable result about the $L_p$ affine surface area is the characterization theorem: roughly speaking, any upper (lower) semi-continuous and affine invariant  (with certain homogeneous degree)  valuation on convex bodies is essentially (up to a multiplicative constant) the $L_p$ affine surface area for some $p>0$  ($-n<p<0$) \cite{LuR, LudR}.  Contributions on the $L_p$ affine surface area include  \cite{lei, Lu91, MW1, MW,   SW4, SWerner, W,  WY} among others. In particular, 
the $L_p$ affine surface area  plays fundamental roles  in the theory of valuations
(see e.g.  \cite{A1, A2, LuR, LudR}), in approximation of convex bodies by
polytopes (see e.g. \cite{Gr2,  LSW, SWerner}) and in the information theory of convex bodies (see e.g.,
\cite{Jenkinson2012, Paouris2010, Werner2012a, Werner2012b}). A full set of affine
isoperimetric inequalities related to the $L_p$ affine surface area has been established  in \cite{Lu96,
WY}. That is, among convex bodies with fixed volume and with centroid at the origin, the $L_p$ affine surface area attains the maximum (minimum) at and only at origin-symmetric ellipsoids for $p>0$ ($-n<p<0$). 

Closely related to the $L_p$ affine surface area, the $L_p$ geominimal surface area for $p\geq 1$  introduced in \cite{Lu96, Petty74} has many nice properties similar to those of the $L_p$ affine surface area, such as affine invariance with the same homogeneous degrees. It is well-known that the $L_p$ geominimal surface for $p\geq 1$ is continuous but the $L_p$ affine surface area is only upper semi-continuous. Hence these two closely related concepts are different from each other. Moreover, the $L_p$ geominimal surface area does not have a ``nice" integral expression similar to that for the $L_p$ affine surface area (which is essential for extending the $L_p$ affine surface area from $p\geq 1$ to all $-n\neq p\in \bbR$). Recently, motivated by an equivalent formula for the $L_p$ affine surface area, the $L_p$ geominimal surface area was successfully extended to all $-n\neq p\in \bbR$ by the first author in \cite{Y2}. When $p=1$, one gets the classical geominimal surface area \cite{Petty74}, which serves as the bridge between the affine geometry, relative geometry and Minkowski geometry (as claimed in \cite{Petty74}). Contributions include \cite{Lu96,  Petty74, Petty1985,  Sch, Y2, zhb} among others. Affine isoperimetric inequalities related to the $L_p$ geominimal surface area can be found in  \cite{Lu96, Petty74, Petty1985, Y2}. That is,  origin-symmetric ellipsoids are the only maximizers (minimizers) of the $L_p$ geominimal surface area  for $p>0$ (for $-n<p<0$), among convex bodies with fixed volume and with centroid at the origin.

A well-studied concept in the literature of convex geometry is the $L_p$ affine surface area for multiple convex bodies, named as the mixed $p$-affine surface area (see e.g., \cite{Lut1987, Lu96, WL1, WY2}). The mixed $p$-affine surface area contains many important functionals on convex bodies as special cases, such as the $L_p$ affine surface area and the dual mixed volume.  From the information theory point of view, the mixed $p$-affine surface area is a very special $f$-divergence of the distributions associated with multiple convex bodies \cite{Werner2012b, WernerYe2013}, which can be used to measure the similarity between multiple convex bodies. Note that the definition of the mixed $p$-affine surface area is also based on a nice integral expression for the $L_p$ affine surface area. Moreover, such a nice integral expression makes it possible to define the general mixed affine surface areas involving nonhomogenous convex and concave functions \cite{Y}.  Alexander-Fenchel type
inequalities  and affine isoperimetric
inequalities for the mixed $p$-affine surface area were established in  \cite{Lu75, Lut1987, Lu96, WY2}.  On the other hand, by using the so-called $p$-Petty body for $p\geq 1$, the second and the third authors propose a way to define the mixed $L_p$ geominimal surface area for multiple convex bodies \cite{zhuzhou} for $p\geq 1$. Related Alexander-Fenchel type
inequalities  and affine isoperimetric
inequalities were also established in \cite{zhuzhou}. 

In this paper, we further extend the mixed $L_p$ geominimal surface area to all $-n\neq p\in \bbR$. Our definitions are motivated from an equivalent formula  for the mixed $p$-affine surface area (see Proposition \ref{equal:formula} and Theorem \ref{equivalent:affine:surface:area-2}), which differ from the definition in \cite{zhuzhou}. Similar idea was used to successfully extend the $L_p$ geominimal surface area to all $-n\neq p\in \bbR$ in \cite{Y2}. Our paper is organized as follows. Section \ref{section:background} is for background and notation. In Section \ref{section:3}, we prove Proposition \ref{equal:formula} and Theorem \ref{equivalent:affine:surface:area-2}. Hence, an equivalent formula  for the mixed $p$-affine surface area is provided.  In particular, the mixed $p$-affine surface area could be rewritten as the dual mixed  volume of the $p$-curvature images of corresponding convex bodies. Section \ref{section:mixed} is dedicated to the definition of our mixed $L_p$ geominimal surface areas for multiple convex bodies for all $-n\neq p\in \bbR$. Properties such as affine invariance are proved. Moreover, relations between mixed $L_p$ geominimal surface areas and the mixed $p$-affine surface area are discussed. Alexander-Fenchel type inequality, affine isoperimetric inequality, the Santal\'{o} style inequality, and the cyclic inequality are proved in Section \ref{section:inequalities}. The $i$-th mixed $L_p$ geominimal surface areas and related isoperimetric inequality are given in Section \ref{section:ith}. We refer readers to  \cite{Gruber2007, Lei1980, Schn} for more background in convex geometry.

\section{Background and Notation}\label{section:background}
We will work on $\bbR^n$ with inner product $\langle \cdot, \cdot\rangle$ and the Euclidean norm $\|\cdot\|$. We use $\ball=\{x\in \bbR^n: \|x\|\leq 1\}$ and
$S^{n-1}=\{x\in \bbR^n: \|x\|= 1\}$ for the unit Euclidean ball
 and the unit sphere in $\mathbb{R}^{n}$, respectively.  For a subset $K\subset \bbR^n$, its Hausdorff content is denoted by $|K|$. In particular, the volume of $\ball$  is
written as $\omega_n=|\ball|$.

A set $L\subset\mathbb{R}^n$ is star-convex about the origin $0$ if for each $x\in L$, the line segment from $0$ to $x$ is contained in $L$. The radial function of $L$, denoted by $\r_L: S^{n-1}\rightarrow [0, \infty)$, is
defined by $\r_L(u)=\max\{\lambda \geq 0:\lambda u\in L\}.$ If
$\rho_L$ is positive and continuous, then $L$ is called a star-convex body
about the origin. The set of all star-convex bodies about the origin is denoted by $\mathcal{S}_0$. We say that $L_1, L_2\in \cS_0$ are dilates of one another if there is a constant $\l>0$, such that $\r_{L_1}(u)=\l \r_{L_2}(u)$ for all $u\in S^{n-1}$. The volume of $L\in \cS_0$ can be calculated by 
\begin{equation*}
|L|=\frac{1}{n}\int_{S^{n-1}}\rho^{n}_{L}(u)d\s(u),
\end{equation*} where $\s$ is the spherical measure on $S^{n-1}$. 

 We say that $K\subset \bbR^n$ is a convex body if $K$ is a compact,
convex subset in $\mathbb{R}^{n}$ with non-empty interior. The set of all convex bodies is written as $\cK$, and its subset $\mathcal{K}_{0}$ denotes the set of convex bodies containing the origin in their interiors. Similarly, we use $\mathcal{K}_{c}$ for the
set of convex bodies with centroid at the origin. Besides the radial function, any convex body can be uniquely determined by its support function. Here, for $K\in\mathcal{K}_0$, its
support function
$h_K: S^{n-1} \rightarrow [0,\infty)$ is
defined by $h_K(u)=\max\{\langle
x, u\rangle: x\in K\}.$ Associated with each $K\in \cK_0$, one can uniquely define its polar body $\polar\in \cK_0$ by $$\polar= \{x\in\mathbb{R}^{n}: \langle x, y\rangle \leq
1,\ \ \  \forall y\in K \}.$$ It is easily verified that
$(K^\circ)^\circ=K$ for all $K\in\mathcal {K}_0$. Moreover,
\begin{equation*}
h_{K^{\circ}}(u) {\rho_{K}(u)}= {1}\ \  \& \ \ \rho_{K^{\circ}}(u){h_{K}(u)}=1, \quad \text{for
all}\quad u\in S^{n-1}.
\end{equation*} A convex body $K\in \cK_0$ is said to have Santal\'{o} point at the
origin, if $\polar$ has centroid at the origin,
i.e., $K \in \cK_s \Leftrightarrow K^\circ \in \cK_c$. Hereafter $\cK_s \subset \cK_0$ denotes the set of all convex bodies with Santal\'{o} point at the origin.

 The $p$-mixed volume, $V_p(K,L)$, of $K, L\in
 \mathcal{K}_0$ for $p\geq1$ was defined in \cite{Lu93} by
 \begin{equation*}
 \frac{n}{p}V_p(K,L)=\lim\limits_{\varepsilon\rightarrow0}
 \frac{|K+_p\varepsilon L|-|K|}{\varepsilon}, 
 \end{equation*} where $K+_p\varepsilon L$ is a convex body with support function defined by 
  \begin{equation*}
 \big(h_{K+_p \varepsilon L}(u)\big)^p= \big(h_K(u)\big)^p+\varepsilon
 \big (h_L(u)\big)^p, \ \ \ \ \forall u\in S^{n-1}.
 \end{equation*} This sum is the well-known Firey
 $p$-sum \cite{Firey}, which generalizes the famous Minkowski sum (i.e., $p=1$). Note that the Minkowski sum $\lambda K+\eta L$ for $K, L\in \cK_0$ and for $\lambda, \eta>0$ is defined by $$h_{\lambda K+\eta L}(u)=\lambda h_K(u)+\eta h_L(u), \ \ \forall u\in S^{n-1}.$$  
 It is well-known that for all $\l_1, \cdots, \l_m>0$ and $K_1, \cdots, K_m\in \cK_0$ with $m\in \mathbb{N}$, one has $$|\l_1K_1+\cdots + \l_m K_m| =\sum_{ i_1,\cdots, i_n=1}^m \l_{i_1} \cdots \l_{i_n} V(K_{i_1}, \cdots, K_{i_n}).$$ The coefficient $V(K_{i_1}, \cdots, K_{i_n})$,  named as the mixed volume of $K_{i_1}, \cdots, K_{i_n}$, is invariant under permutations of $K_{i_1}, \cdots, K_{i_n}$. The classical Alexander-Fenchel inequality for the mixed volume  (see \cite{LSW, Schn}) asserts that for all $m\in \mathbb{N}$ such that $1\leq m\leq n$,
 \begin{equation*}
 \prod_{i=0}^{m-1}V(K_{1},\cdots, K_{n-m},
 \underbrace{K_{n-i},\cdots, K_{n-i}}_{m}) \leq V(K_{1},\cdots,
 K_{n})^{m}.
 \end{equation*} In particular, if $m=n$, one has the Minkowski
 inequality  
 \be \label{Minkowski inequality-mixed volume} V(K_{1},\cdots,
 K_{n})^{n}\geq |K_{1}|\cdots |K_{n}|. \ee

  It was proved in \cite{Lu93} that for
each $K\in\mathcal{K}_0$, there is a positive Borel measure
$S_{p}(K,\cdot)$ on $S^{n-1}$ such that, for each $L\in\mathcal{K}_0$ and $p\geq 1$,
 \be
\nonumber V_{p}(K,
L)=\frac{1}{n}\int_{S^{n-1}}h_L(u)^{p}dS_{p}(K, u).\ee
 Moreover, the measure $S_p(K,\cdot)$ for $p\geq 1$ has the following form \be\nonumber
 {dS_p(K,\cdot)}=h_K(\cdot)^{1-p}{dS(K,\cdot)}, \ee where  the measure $S(K,\cdot)$ is
just the classical surface area measure of $K$ (see \cite{Ale1937-1, Fenchel1938}). We write $K\in \cF_0$ if $K\in \cK_0$ has a curvature function, namely, the measure $S(K, \cdot)$ is absolutely continuous with respect to the spherical measure $\s$. Hence, there is a function $f_K: S^{n-1}\rightarrow \bbR$, the {\it curvature function} of $K$, such that, $$dS(K,u)=f_K(u)d\s(u).$$ The {\it $L_p$ curvature function} for $K\in \cF_0$ and $p\geq 1$, denoted by $f_p(K, u)$ (see \cite{Lu96}) then takes the form  $$f_p(K, u)=h_K(u)^{1-p}f_K(u).$$ We write $\cF_c=\cF_0\cap \cK_c$ and $\cF_s=\cF_0\cap \cK_s$ for convex bodies in $\cF_0$ with centroid and the Santal\'{o} point at the origin respectively. The set of all convex bodies in $\cF_0$ with continuous positive curvature function $f_K(\cdot)$  on $S^{n-1}$ is denoted by $\cF_0^+$. 

The mixed $p$-volume for $p\geq 1$ was formally extended to  $p<1$ in \cite{Y2}. Hereafter, the $p$-surface area measure of $K\in \cK_0$ is
$$dS_p(K,u)=h_K(u)^{1-p}dS(K, u), \ \ \ p\in \bbR,$$ and the $p$-mixed
volume of $K, Q\in \cK_0$ is
 \be\label{p-volume-integral-1-----1} V_p(K, Q)=\frac{1}{n}\int_{S^{n-1}}h_Q(u)^p dS_p(K,u), \ \ \ p\in \bbR. \ee
 For $K\in \cK_0$ and $L\in \cS_0$, we let $V_p(K,L^\circ)$ be 
\be \nonumber V_p(K,L^\circ)=\frac{1}{n}\int_{S^{n-1}} \r_L(u)^{-p}dS_p(K, u), \ \ \ p\in \bbR.
 \ee This formula coincides with formula
    (\ref{p-volume-integral-1-----1}) if $L\in \cK_0$. When $K\in \cF_0$, one can define the $L_p$ curvature function (and denoted by $f_p(K, \cdot)$) as
\be\nonumber 
f_p(K, u)=h_K(u)^{1-p}f_K(u),  \ \ \ p\in \bbR,
\ee
 and hence the $p$-surface area measure can be formulated by
\be\nonumber dS_p(K,u)=f_{p}(K,u)d\sigma(u), \ \ \ p\in \bbR.
\ee

 Now we can define the $L_p$ geominimal surface area of $K\in \cK_0$ as follows.  See \cite{Petty74} for $p=1$, \cite{Lu96} for $p>1$ and \cite{Y2} for all $-n\neq p\in \bbR$. 
  
  \bd \label{p-geominimal}
  Let $K\in \cK_0$ be a convex body with  the origin in its interior.   \\
  (i). For $p=0$,  we let $\widetilde{G}_p(K)=n |K|$. For $p>0$,  the
  $L_p$  geominimal surface area of $K$ is defined by
   \begin{equation*} \widetilde{G}_p(K)=\inf_{L\in \cK_0 }
   \left\{n V_p(K,  L) ^{\frac{n}{n+p}}\
  |L ^\circ|^{\frac{p}{n+p}}\right\} .
  \end{equation*}
  (ii). For $-n\neq p<0$, the $L_p$ geominimal surface area of $K$ is
  defined by
  \begin{equation*}\widetilde{G}_p(K)=\sup_{L \in \cK_0}
  \left\{nV_p(K,  L) ^{\frac{n}{n+p}}\ |L
  ^\circ|^{\frac{p}{n+p}}\right\} .
  \end{equation*}
  \ed
  
We let $\cL=(L_1, \cdots, L_n)$ be a vector with each $L_i\subset \bbR^n$, and $\cL\in \cS_0^n$ means  each $L_i\in \cS_0$. We use  $\cL^\circ$ for $(L_1^\circ, \cdots, L_n^\circ)$. Similarly, we let $\bK=(K_1, \cdots, K_n)$ and $\bK\in \cF_0^n$ means each $K_i\in \cF_0$. For all $\{K_i\}_{i=1}^n\subset \cF_0$, $\{Q_i\}_{i=1}^n\subset \cK_0$, and  $p\in \bbR$, we define \be \label{p-volume:convex} V_p(\bK; \bQ) =\frac{1}{n}
\int_{S^{n-1}} \prod _{i=1}^n \left[h_{Q_i} (u)^{p} f_p(K_i,
u)\right]^{\frac{1}{n}}\,d\s(u).  \ee When all $K_i$ coincide with $K$ and all $Q_i$ coincide with $Q$, one can easily get \be V_p(\bK; \bQ)=V_p(K, Q).\nonumber \ee  When  $L_1,
\cdots, L_n\in \cS_0$ and $p\in \bbR$, we use the following variation formula \be \label{p-volume} V_p(\bK; \cL^\circ) =\frac{1}{n}
\int_{S^{n-1}} \prod _{i=1}^n \left[\r _{L_i} (u)^{-p} f_p(K_i,
u)\right]^{\frac{1}{n}}\,d\s(u),  \ee which is consistent with  formula  (\ref{p-volume:convex}) when $\cL\in \cK_0^n$.  When all $K_i$ coincide with $K$ and all $L_i$ coincide with $L$, one gets \be V_p(\bK; \cL^\circ)=V_p(K, L^\circ).\label{coincide:mixed}\ee
We use $\widetilde{V}(\cL)$ to denote the dual mixed volume  
of $L_1, \cdots, L_n\in \cS_0$ (see\cite{Lu75}). That is, 
$$\widetilde{V}(\cL)=\widetilde{V}(L_1, \cdots, L_n) =\frac{1}{n}\int_{S^{n-1}}
\prod_{i=1}^n \r_{L_i}(u) \,d\s(u).$$ When $L_1=L_2=\cdots=L_n=L$, one has $\widetilde{V}(\cL)=|L|$. It is easy to get the following
inequality for the dual mixed volume: 
\begin{equation}\label{dual mixed volume}
[\widetilde{V}(\cL)]^n=[\widetilde{V}(L_{1},\cdots, L_{n})]^{n}\leq |L_{1}|\cdots |L_{n}|,
\end{equation}
with equality if and only if  $L_{i}$ ($1\, \le i \, \le n$)  are
dilates of each other.

The group of nonsingular linear
transformations is denoted by $GL(n)$. Its subset $SL(n)$ refers to the group of special linear transformations. For $\phi\in GL(n)$, the absolute
value of the determinant, the transpose and the inverse of $\phi$ will be denoted by $|det (\phi)|$, $\phi^{t}$ and $\phi^{-1}$ respectively. For $\bK=(K_1, \cdots, K_n)\in \cK_0^n$, we let $\phi \bK=(\phi K_1, \cdots, \phi K_n)$. An origin-symmetric ellipsoid $\mathcal{E} \in \cK_0$
can be obtained from $\ball$ under some $\phi\in GL(n)$, that is,  $\mathcal{E} = \phi
\ball$ for some $\phi\in GL(n)$.  Note that $(\phi L)^\circ=(\phi ^t)^{-1} L^\circ$ for $\phi\in GL(n)$ and $L\in \cS_0$, then \be \widetilde{V}((\phi \cL)^\circ)=|det (\phi)|^{-1} \widetilde{V}(\cL^\circ). \label{affine:dual:mixed:volume}\ee

\section{The mixed $p$-affine surface area, another view}\label{section:3}
In this section, we will prove some alternative formulas for the mixed $p$-affine surface area. Our definition for the mixed $L_p$ geominimal surface areas for multiple convex bodies is motivated by these alternative formulas.

\bp\label{equal:formula} Let $K_1, \cdots, K_n\in \cF_0$.  
\\ (i).
For $p\geq 0$,  one has
\begin{eqnarray*}  \inf_{ \cL\in \cS_0^n} \left\{n  V_p(\bK;  \cL^\circ) ^{\frac{n}{n+p}}\ \widetilde{V}(\cL)^{\frac{p}{n+p}}\right\} &=& \inf_{  \cL\in \cS_0^n } \bigg\{n  V_p(\bK;  \cL^\circ) ^{\frac{n}{n+p}}\ \prod_{i=1}^n
|L_i|^{\frac{p}{(n+p)n}}\bigg\} \\ &=& \inf_{ L  \in \cS_0}
 \left\{nV_p(\bK;  L ^\circ, \cdots, L ^\circ) ^{\frac{n}{n+p}}\
|L|^{\frac{p}{n+p}}\right\} .
\end{eqnarray*}
(ii). For $-n<p<0$, one has
\begin{eqnarray*} \sup_{  \cL\in \cS_0^n} \left\{n  V_p(\bK;  \cL^\circ) ^{\frac{n}{n+p}}\
\widetilde{V}(\cL)^{\frac{p}{n+p}}\right\} &=&\sup_{
 \cL\in \cS_0^n } \bigg\{n  V_p(\bK;
\cL^\circ) ^{\frac{n}{n+p}}\ \prod_{i=1}^n
|L_i|^{\frac{p}{(n+p)n}}\bigg\}\\&=&\sup_{ L
\in \cS_0  } \left\{n  V_p(\bK;  L ^\circ, \cdots, L
^\circ) ^{\frac{n}{n+p}}\ |L|^{\frac{p}{n+p}}\right\}.\end{eqnarray*}
(iii). 
For $
p<-n$, one has
 \begin{eqnarray*} \sup_{
 \cL\in \cS_0^n } \left\{n V_p(\bK;
\cL^\circ) ^{\frac{n}{n+p}}\ \widetilde{V}(\cL)^{\frac{p}{n+p}}\right\} =\sup_{ L \in \cS_0  }
\left\{n V_p(\bK;  L ^\circ, \cdots, L ^\circ)
^{\frac{n}{n+p}}\ |L|^{\frac{p}{n+p}}\right\}.
\end{eqnarray*}\ep

\vskip 2mm \noindent {\bf Proof.} First, notice that there is a sequence $\{\cL_j\}_{j=1}^{\infty}\subset\cS_0^n$, s.t., \begin{eqnarray} \inf_{ \cL\in \cS_0^n} \left\{n  V_p(\bK;  \cL^\circ) ^{\frac{n}{n+p}}\ \widetilde{V}(\cL)^{\frac{p}{n+p}}\right\} &=&\lim_{j\rightarrow\infty}  \left\{n  V_p(\bK;  \cL_j^\circ) ^{\frac{n}{n+p}}\ \widetilde{V}(\cL_j)^{\frac{p}{n+p}}\right\}, \ \ p>0; \label{case:p>0} \\ \sup_{ \cL\in \cS_0^n} \left\{n  V_p(\bK;  \cL^\circ) ^{\frac{n}{n+p}}\ \widetilde{V}(\cL)^{\frac{p}{n+p}}\right\} \!\!\!&=&\lim_{j\rightarrow\infty}  \left\{n  V_p(\bK;  \cL_j^\circ) ^{\frac{n}{n+p}}\ \widetilde{V}(\cL_j)^{\frac{p}{n+p}}\right\}, \ \ -n\neq p<0. \ \ \ \ \ \ \ \ \label{case:p<0}  \end{eqnarray}  For each $\cL_j\!=\!(L_{1j}, \cdots, L_{nj})\!\in\! \cS_0^n$, one defines $L^j\in \cS_0$  by $\r_{L^j}^n(u)=\prod_{i=1}^n \r_{L_{ij}}(u).$ Hence,  for all $j$, \be\widetilde{V}(\cL_j)=\frac{1}{n}\int_{S^{n-1}}
\prod_{i=1}^n \r_{L_{ij}}\,d\s(u)=\frac{1}{n}\int_{S^{n-1}}
\r_{L^j}(u)^n\,d\s(u)=|L^j|. \label{dual:mixed:surface-1}\ee Moreover, by formula (\ref{p-volume}), one has, \begin{eqnarray}
V_p(\bK; \cL_j^\circ)  &=&\frac{1}{n}
\int_{S^{n-1}} \prod _{i=1}^n  [\r _{L_{ij}} (u)^{-p}]^{\frac{1}{n}} \prod _{i=1}^n [f_p(K_i,
u) ]^{\frac{1}{n}}\,d\s(u)\nonumber\\&=& \frac{1}{n}
\int_{S^{n-1}}  \r _{L^j } (u)^{-p} \prod _{i=1}^n [f_p(K_i,
u) ]^{\frac{1}{n}}\,d\s(u)= V_p(\bK;  (L^j)^\circ, \cdots, (L^j)^\circ).\label{dual:mixed:surface-2}
\end{eqnarray} 
(i). The case of $p=0$ is clear.  Let $p\in (0, \infty)$, then \begin{eqnarray*}  \inf_{ \cL\in \cS_0^n} \left\{n  V_p(\bK;  \cL^\circ) ^{\frac{n}{n+p}}\ \widetilde{V}(\cL)^{\frac{p}{n+p}}\right\} &\leq& \inf_{  \cL\in \cS_0^n } \bigg\{n  V_p(\bK;  \cL^\circ) ^{\frac{n}{n+p}}\ \prod_{i=1}^n
|L_i|^{\frac{p}{(n+p)n}}\bigg\} \\ &\leq& \inf_{ L  \in \cS_0}
 \left\{nV_p(\bK;  L ^\circ, \cdots, L ^\circ) ^{\frac{n}{n+p}}\
|L|^{\frac{p}{n+p}}\right\},
\end{eqnarray*}where the first inequality follows from inequality (\ref{dual mixed volume}) and the
second inequality is due to  $(L,
\cdots, L)\in \cS_0^n$. On the other hand, formulas (\ref{case:p>0}), (\ref{dual:mixed:surface-1}) and (\ref{dual:mixed:surface-2})  imply that \begin{eqnarray*} \inf_{ \cL\in \cS_0^n} \left\{n  V_p(\bK;  \cL^\circ) ^{\frac{n}{n+p}}\ \widetilde{V}(\cL)^{\frac{p}{n+p}}\right\} \!\!\!&=&\!\!\!\lim_{j\rightarrow\infty}  \left\{n  V_p(\bK;  \cL_j^\circ) ^{\frac{n}{n+p}}\ \widetilde{V}(\cL_j)^{\frac{p}{n+p}}\right\}\\\!\!\!&=&\!\!\!\lim_{j\rightarrow\infty}  \left\{n  V_p(\bK;  (L^j)^\circ, \cdots, (L^j)^\circ) ^{\frac{n}{n+p}}\ |L^j|^{\frac{p}{n+p}}\right\}\\&\geq&\inf_{ L  \in \cS_0}
 \left\{nV_p(\bK;  L ^\circ, \cdots, L ^\circ) ^{\frac{n}{n+p}}\
|L|^{\frac{p}{n+p}}\right\}.\end{eqnarray*}
(ii).  Let $p\in (-n, 0)$, then $\frac{p}{n+p}<0$. By formulas (\ref{case:p<0}), (\ref{dual:mixed:surface-1}) and (\ref{dual:mixed:surface-2}),  inequality (\ref{dual mixed volume}), and $(L,
\cdots, L)\in \cS_0^n$, one has \begin{eqnarray*}  \sup_{ \cL\in \cS_0^n} \left\{n  V_p(\bK;  \cL^\circ) ^{\frac{n}{n+p}}\ \widetilde{V}(\cL)^{\frac{p}{n+p}}\right\} &\geq& \sup_{  \cL\in \cS_0^n } \bigg\{n  V_p(\bK;  \cL^\circ) ^{\frac{n}{n+p}}\ \prod_{i=1}^n
|L_i|^{\frac{p}{(n+p)n}}\bigg\} \\ &\geq& \sup_{ L  \in \cS_0}
 \left\{nV_p(\bK;  L ^\circ, \cdots, L ^\circ) ^{\frac{n}{n+p}}\
|L|^{\frac{p}{n+p}}\right\}\\ &\geq & \!\!\!\lim_{j\rightarrow\infty}  \left\{n  V_p(\bK;  (L^j)^\circ, \cdots, (L^j)^\circ) ^{\frac{n}{n+p}}\ |L^j|^{\frac{p}{n+p}}\right\}  \\ &=&  \sup_{ \cL\in \cS_0^n} \left\{n  V_p(\bK;  \cL^\circ) ^{\frac{n}{n+p}}\ \widetilde{V}(\cL)^{\frac{p}{n+p}}\right\} .\end{eqnarray*}
(iii). Let $p<-n$. Due to  $(L,
\cdots, L)\in \cS_0^n$, and  formulas (\ref{case:p<0}), (\ref{dual:mixed:surface-1}) and (\ref{dual:mixed:surface-2}),  one has  \begin{eqnarray*}  \sup_{ \cL\in \cS_0^n} \left\{n  V_p(\bK;  \cL^\circ) ^{\frac{n}{n+p}}\ \widetilde{V}(\cL)^{\frac{p}{n+p}}\right\} &\geq& \sup_{ L  \in \cS_0}
 \left\{nV_p(\bK;  L ^\circ, \cdots, L ^\circ) ^{\frac{n}{n+p}}\
|L|^{\frac{p}{n+p}}\right\}\\ &\geq & \lim_{j\rightarrow\infty}  \left\{n  V_p(\bK;  (L^j)^\circ, \cdots, (L^j)^\circ) ^{\frac{n}{n+p}}\ |L^j|^{\frac{p}{n+p}}\right\}\\ &=& \sup_{ \cL\in \cS_0^n} \left\{n  V_p(\bK;  \cL^\circ) ^{\frac{n}{n+p}}\ \widetilde{V}(\cL)^{\frac{p}{n+p}}\right\}. \end{eqnarray*}

\noindent {\bf Remark.} For $\bK\in \cF_0^n$ and for $p<-n$, due to $(L, \cdots, L)\in \cS_0^n$, one can prove that \begin{eqnarray} \sup_{
 \cL\in \cS_0^n } \bigg\{n  V_p(\bK;
\cL^\circ) ^{\frac{n}{n+p}}\ \prod_{i=1}^n
|L_i|^{\frac{p}{(n+p)n}}\bigg\}&\geq &\sup_{ L
\in \cS_0  } \left\{n  V_p(\bK;  L ^\circ, \cdots, L
^\circ) ^{\frac{n}{n+p}}\ |L|^{\frac{p}{n+p}}\right\}\nonumber \\ &=& \sup_{  \cL\in \cS_0^n} \left\{n  V_p(\bK;  \cL^\circ) ^{\frac{n}{n+p}}\
\widetilde{V}(\cL)^{\frac{p}{n+p}}\right\} . \label{not equal:p<-n}\end{eqnarray} However, $\exists \bK\in \cF_0^n$, s.t. ``$>$" holds (see the remark after Definition \ref{equivalent:mixed affine:surface:area-1} for more details). 

In literature, the mixed $p$-affine surface area of convex bodies $K_1, \cdots, K_n\in \cF_0^+$, denoted by $as_{p}(K_1, \cdots, K_n)$, was defined by \cite{Lu96, WY2}  \be
as_{p}(K_1, \cdots, K_n)=\int_{S^{n-1}} {\prod _{i=1}^n [f_p({K_i},u)]^{\frac{1}{n+p}}} d\sigma(u). \label{Lp-affine-surface-area-integral} \ee 
When all $K_i$ coincide with $K\in \cF_0^+$, one gets, by Theorem 3.1 in \cite{Y2},  $$as_p(K_1, \cdots, K_n)=as_p(K)=\left\{ \begin{array}{cl}
\inf_{ L  \in \cS_0}
 \left\{n  V_p(K,  L ^\circ) ^{\frac{n}{n+p}}\
|L|^{\frac{p}{n+p}}\right\},& p\geq 0; \\ 
\sup_{ L  \in \cS_0}
 \left\{n V_p(K,  L ^\circ) ^{\frac{n}{n+p}}\
|L|^{\frac{p}{n+p}}\right\}, & -n \neq p< 0. 
\end{array}  \right.$$
 
\bt\label{equivalent:affine:surface:area-2}  Let $K_1, \cdots, K_n\in \cF_0^+$.   \\
(i). For $p\geq 0$,  one has
\begin{eqnarray*} as_p(K_1, \cdots, K_n) &=&\inf_{ L  \in \cS_0}
 \left\{n  V_p(\bK;  L ^\circ, \cdots, L ^\circ) ^{\frac{n}{n+p}}\
|L|^{\frac{p}{n+p}}\right\} .
\end{eqnarray*}
(ii). For $-n\neq p<0$, one has
\begin{eqnarray*} as_p(K_1, \cdots,
K_n)&=&\sup_{ L
\in \cS_0  } \left\{n  V_p(\bK;  L ^\circ, \cdots, L
^\circ) ^{\frac{n}{n+p}}\ |L|^{\frac{p}{n+p}}\right\}.\end{eqnarray*}
\et

\noindent \textbf{Proof.}  First, notice that for all $-n\neq p\in \bbR$, \begin{eqnarray} as_p(K_1,
\cdots, K_n) &=& n V_p(\bK;  L_{0} ^\circ, \cdots, L_{0} ^\circ) ^{\frac{n}{n+p}}\
|L_{0}|^{\frac{p}{n+p}}, \label{mixed p curvature-1-L}\end{eqnarray}  where   $L_0\in \cS_0$ is defined by 
$$\r_{L_0}^n=\bigg(\prod _{i=1}^n
f_p({K_i}, u)\bigg)^{\frac{1}{n+p}}>0, \ \ \forall u\in S^{n-1}. 
$$ 
(i). Clearly it holds for $p=0$. Let $p\in (0,\infty)$, then $\frac{n}{n+p}\in (0,1)$. Employing H\"{o}lder inequality (see \cite{HLP}) to formula (\ref{Lp-affine-surface-area-integral}), one has, for all $ L \in \cS_0$,
\begin{eqnarray*}
as_p(K_1, \cdots, K_n)&=&\int_{S^{n-1}}{\bigg(\prod _{i=1}^n
[\r_{L}(u)^{-p}f_p({K_i}, u)]
^{\frac{1}{n}}\bigg)^{\frac{n}{n+p}}}\bigg( \r_{L}^n(u)\bigg)^{\frac{p}{n+p}} d\sigma(u)\nonumber  \\
&\leq&\bigg(\!\int_{S^{n-1}}{\prod _{i=1}^n
[\r_{L}(u)^{-p}f_p({K_i}, u)]
^{\frac{1}{n}}d\sigma(u)}\!\bigg)^{\frac{n}{n+p}}\bigg(\!\int_{S^{n-1}}{\r_{L}^n(u)}d\sigma(u)\!\bigg)^{\frac{p}{n+p}} \nonumber \\
&=& n V_p(\bK;  L ^\circ, \cdots, L ^\circ)
^{\frac{n}{n+p}}\ |L|^{\frac{p}{n+p}}. \end{eqnarray*} Taking infimum over $L \in \cS_0$ and together with formula (\ref{mixed p curvature-1-L}), one gets 
 $$as_p(K_1, \cdots, K_n)\leq\inf_{L\in\cS_0} \left\{n V_p(\bK;  L ^\circ, \cdots, L ^\circ)
^{\frac{n}{n+p}}\ |L|^{\frac{p}{n+p}}\right\}\leq as_p(K_1, \cdots, K_n).$$

\vskip 2mm \noindent (ii). Note that $-n\neq p<0$
 implies $\frac{n}{n+p}>1$ or $\frac{n}{n+p}<0$. Employing H\"{o}lder inequality (see \cite{HLP}) to formula (\ref{Lp-affine-surface-area-integral}), one has, for all $L \in \cS_0$,
\begin{eqnarray*}
as_p(K_1, \cdots, K_n)&=&\int_{S^{n-1}}{\bigg(\prod _{i=1}^n
[\r_{L}(u)^{-p}f_p({K_i}, u)]
^{\frac{1}{n}}\bigg)^{\frac{n}{n+p}}}\bigg( \r_{L}^n(u)\bigg)^{\frac{p}{n+p}} d\sigma(u)\nonumber  \\
&\geq&\bigg(\!\int_{S^{n-1}}{\prod _{i=1}^n
[\r_{L}(u)^{-p}f_p({K_i}, u)]
^{\frac{1}{n}}d\sigma(u)}\!\bigg)^{\frac{n}{n+p}}\bigg(\!\int_{S^{n-1}}{\r_{L}^n(u)}d\sigma(u)\!\bigg)^{\frac{p}{n+p}} \nonumber \\
&=& n V_p(\bK;  L ^\circ, \cdots, L ^\circ)
^{\frac{n}{n+p}}\ |L|^{\frac{p}{n+p}}.
\end{eqnarray*} Taking supremum over $L \in \cS_0$ and together with formula (\ref{mixed p curvature-1-L}), one gets 
 $$as_p(K_1, \cdots, K_n)\geq\sup_{L\in\cS_0} \left\{n V_p(\bK;  L ^\circ, \cdots, L ^\circ)
^{\frac{n}{n+p}}\ |L|^{\frac{p}{n+p}}\right\}\geq as_p(K_1, \cdots, K_n).$$

Motivated by Proposition \ref{equal:formula} and Theorem
\ref{equivalent:affine:surface:area-2}, the mixed $p$-affine surface area for $K_1, \cdots, K_n\in \cF_0$ may be defined as follows.

\bd\label{equivalent:mixed affine:surface:area-1} Let $K_1, \cdots,
K_n\in \cF_0$. \\
(i). For $p=0$,  we let $$as_0 (K_1, \cdots, K_n)=\int_{S^{n-1}}
{\prod _{i=1}^n [f_{K_i}(u)h_{K_i}(u)]^{\frac{1}{n}}} d\sigma(u).$$ For
$p>0$, we define  $as_p (K_1, \cdots, K_n)$ by
 \begin{equation*} as_p (K_1, \cdots, K_n)=\inf_{ L\in \cS_0  }
 \left\{n  V_p (\bK;  L^\circ, \cdots, L ^\circ) ^{\frac{n}{n+p}}\
 |L|^{\frac{p}{n+p}}\right\} .
\end{equation*}
(ii). For $-n\neq p<0$, we define  $as_p(K_1, \cdots, K_n)$  by
\begin{equation*} as_p (K_1, \cdots, K_n)=\sup_{  L \in \cS_0  }
\left\{n V_p(\bK;  L ^\circ, \cdots, L ^\circ)
^{\frac{n}{n+p}}\ |L|^{\frac{p}{n+p}}\right\} .
\end{equation*}
Moreover, for $p<-n$, one can define another mixed $p$-affine surface area for $\bK\in \cF_0^n$ as 
\begin{equation*} as_p^{(1)} (K_1, \cdots, K_n)=\sup_{ \cL \in \cS_0^n }
\bigg\{n  V_p(\bK;  \cL ^\circ)
^{\frac{n}{n+p}}\ \prod_{i=1}^n |L_i|^{\frac{p}{(n+p)n}}\bigg\}.
\end{equation*}
\ed
\vskip 2mm \noindent {\bf Remark.} For $\bK\in \cF_0^n$ and for $p<-n$, inequality (\ref{not equal:p<-n}) implies   \begin{eqnarray*} as_p (K_1, \cdots, K_n) &\leq& as_p^{(1)} (K_1, \cdots, K_n). \end{eqnarray*}  In general, one {\it cannot} expect equality in the above inequality. To see this, suppose that $K_1, \cdots, K_n\in \cF_0^+$. Proposition 2.1 in \cite{WY2} and Theorem \ref{equivalent:affine:surface:area-2} imply $$as_p (K_1, \cdots, K_n)  \leq   as_p(K_1)  ^{\frac{1}{n}}\cdots as_p(K_n)  ^{\frac{1}{n}},  $$ with equality if and only if there are constants $\lambda_1, \lambda_2, \cdots, \lambda_n>0$ such that $\lambda_i f_p(K_i, u)=\lambda_j  f_p(K_j, u)$ for all $u\in S^{n-1}$. Now suppose that there is {\it no constant} $\lambda>0$ satisfying  $f_p(K_1, u)= \l f_p(K_2, u)$ almost everywhere with respect to the spherical measure $\s$, then $$as_p (K_1, \cdots, K_n) < as_p(K_1)  ^{\frac{1}{n}}\cdots as_p(K_n)  ^{\frac{1}{n}}.  $$
 On the other hand, by H\"{o}lder inequality (see \cite{HLP}), one gets $$V_p(\bK;  \cL^\circ)\leq  \prod_{i=1}^n \big[V_p(K_i, L_i^\circ)\big]^{\frac{1}{n}}. $$ Note that $p<-n$ implies $\frac{n}{n+p}<0$. Therefore, \begin{eqnarray*} as_p^{(1)} (K_1, \cdots, K_n)&=&\sup_{ \cL \in \cS_0^n }
 \bigg\{n  V_p(\bK;  \cL ^\circ)
 ^{\frac{n}{n+p}}\ \prod_{i=1}^n |L_i|^{\frac{p}{(n+p)n}}\bigg\} \\ 
 &\geq&  \prod_{i=1}^n \sup_{
 L_i \in \cS_0 }    \big[n V_p(K_i, L_i^\circ) ^{\frac{n}{n+p}}  
|L_i|^{\frac{p}{n+p}}\big] ^{\frac{1}{n}}\\ &=& as_p(K_1)  ^{\frac{1}{n}}\cdots as_p(K_n)  ^{\frac{1}{n}} > as_p (K_1, \cdots, K_n). \end{eqnarray*} 

 The $p$-curvature image of $K$ for $-n\neq p\in \bbR$ \cite{Lu96, Y2}  is denoted by $\L_p K$ and defined by $$f_p(K, u)= \frac{\o _n}{ |\L _p K|}\r_{\L _p
K}(u)^{n+p}, \ \ \forall u \in S^{n-1}.$$
 
\bp \label{affine:p:curvature} For $K_1, \cdots, K_n\in
\cF_0^+$ and for $-n\neq p\in \bbR$,
\begin{equation*}
as_{p} (K_{1},\cdots,  K_{n}) ^ {n+p}=\frac{n^{n+p}\omega^n_n}{|\L_p
K_1| \cdots |\L_p K_n|} \widetilde{V}(\L_p K_1, \cdots, \L_p
K_n)^{n+p}.
\end{equation*}
\ep
\noindent{\bf Proof.} Note that formula (\ref{mixed p curvature-1-L}) still holds if one changes $L_0$ to $\lambda L_0$ for any $\lambda>0$. In particular, we choose $L_0$ to be  
\begin{eqnarray*}
\big[\r_ {L_0}(u)\big]^n=\r_{\L _p K_1} (u) \cdots \r_{\L _p K_n}(u)=\bigg(\prod_{i=1}^n \frac{|\Lambda_p K_i|}{\o_n}\cdot \prod _{i=1}^n
f_p({K_i}, u)\bigg)^{\frac{1}{n+p}}, \ \  \forall u\in S^{n-1}.
\end{eqnarray*}
Combining with formula (\ref{Lp-affine-surface-area-integral}), one has  \begin{eqnarray*} n\widetilde{V}(\L _p K_1, \cdots, \L _p K_n)  &=&\int_{S^{n-1}} \r_{\L _p K_1} (u) \cdots \r_{\L _p K_n}(u) \,d\s(u) \\ &=& \int_{S^{n-1}} \bigg(\prod_{i=1}^n \frac{|\Lambda_p K_i|}{\o_n}\cdot \prod _{i=1}^n
f_p({K_i}, u)\bigg)^{\frac{1}{n+p}}\,d\s(u) \\  &=& \bigg(\prod_{i=1}^n \frac{|\Lambda_p K_i|}{\o_n}\bigg)^{\frac{1}{n+p}}  \int_{S^{n-1}} \prod _{i=1}^n
\big[f_p({K_i}, u)\big]^{\frac{1}{n+p}}\,d\s(u) \\&=& \bigg(\prod_{i=1}^n \frac{|\Lambda_p K_i|}{\o_n}\bigg)^{\frac{1}{n+p}}  as_p(K_1, \cdots, K_n). \end{eqnarray*}  This concludes the proof of Proposition \ref{affine:p:curvature}.

\section{Mixed $L_p$ geominimal surface areas}\label{section:mixed}

We now introduce several mixed $L_p$ geominimal
surface areas for all $-n\neq p\in \bbR$.

\bd\label{equivalent:mixed geominimal:surface:area-1}
Let $K_1,\cdots,K_n\in \cF_0$.   \\
(i). For $p=0$,  we let $$G_0^{(\a)} (K_1, \cdots, K_n)=\int_{S^{n-1}}
{\prod _{i=1}^n [f_{K_i}(u)h_{K_i}(u)]^{\frac{1}{n}}} d\sigma(u), \ \ \ \a=1, 2, 3. $$ For
$p>0$,  we define the mixed $L_p$ geominimal surface areas for $\bK \in \cF_0^n$ as:   \begin{eqnarray*} G_p^{(1)}(K_1, \cdots, K_n)&=&\inf_{  L\in \cK_0  }
 \left\{n V_p(\bK;  L, \cdots, L) ^{\frac{n}{n+p}}\
 |L ^\circ|^{\frac{p}{n+p}}\right\}, \\
 G_p^{(2)}(K_1, \cdots, K_n)&=&\inf_{ \cL\in \cK_0^n }
 \bigg\{n V_p(\bK;  \cL) ^{\frac{n}{n+p}}\
\prod_{i=1}^n |L_i ^\circ|^{\frac{p}{(n+p)n}}\bigg\},\\
G_p^{(3)}(K_1, \cdots, K_n)&=&\inf_{ \cL\in \cK_0^n} \left\{n  V_p(\bK;  \cL) ^{\frac{n}{n+p}}\ \widetilde{V}(\cL^\circ)^{\frac{p}{n+p}}\right\} .
\end{eqnarray*}
 
\noindent (ii).\ For $-n\neq p<0$, we define the mixed $L_p$ geominimal surface areas for $\bK \in \cF_0^n$ as:   \begin{eqnarray*} G_p^{(1)}(K_1, \cdots, K_n)&=&\sup_{  L\in \cK_0  }
 \left\{n V_p(\bK;  L, \cdots, L) ^{\frac{n}{n+p}}\
 |L ^\circ|^{\frac{p}{n+p}}\right\}, \\
 G_p^{(2)}(K_1, \cdots, K_n)&=&\sup_{ \cL\in \cK_0^n }
 \bigg\{n V_p(\bK;  \cL) ^{\frac{n}{n+p}}\
\prod_{i=1}^n |L_i ^\circ|^{\frac{p}{(n+p)n}}\bigg\},\\
G_p^{(3)}(K_1, \cdots, K_n)&=&\sup_{ \cL\in \cK_0^n} \left\{n  V_p(\bK;  \cL) ^{\frac{n}{n+p}}\ \widetilde{V}(\cL^\circ)^{\frac{p}{n+p}}\right\} .
\end{eqnarray*}\ed

  We now prove that the mixed $L_{p}$ geominimal surface
areas defined in Definition  \ref{equivalent:mixed
geominimal:surface:area-1}  are all affine invariant.
 \bp \label{proposition-1} If
$K_{1}, \cdots, K_{n}\in\cF_0$ and $\phi\in GL(n)$, then for all $\a=1, 2, 3$, and for all  $-n\neq p\in \bbR$, one has
\begin{equation*}
G_{p}^{(\a)} (\phi K_{1},\cdots, \phi K_{n})=|det
(\phi)|^{\frac{n-p}{n+p}}G_{p}^{(\a)}(K_{1},\cdots, K_{n}).
\end{equation*} 
In particular, $G_{p} ^{(\a)}(K_{1},\cdots,
K_{n})$ for all $\a=1,2, 3$ are affine invariant, namely if $\phi\in SL(n)$, then
$
G_{p}^{(\a)} (\phi K_{1},\cdots, \phi K_{n})=G_{p} ^{(\a)} (K_{1},\cdots,
K_{n}).$
\ep

\noindent \textbf{Proof.}  Let $\phi\in GL(n)$ and $v=\frac{\phi^{-t}(u)}{\|\phi^{-t}(u)\|}$ for any $u \in S^{n-1}$. Then for all $K\in \cK_0$ \be  
 h_{K}(u)=\max_{x\in K} \langle x,
u\rangle=\max_{x\in K}\langle \phi x, \phi^{-t}(u)\rangle=\max_{x\in K}\|\phi^{-t}(u)\|\langle
\phi x, v\rangle=\|\phi^{-t}(u)\| h_{\phi K}(v). \nonumber \ee    Hence, for all $u\in S^{n-1}$, \be \frac{h_{\phi Q_i}(v)}{h_{\phi K_i}(v)}=\frac{h_{Q_i}(u)}{h_{K_i}(u)}. \label{support:ratio} \ee On the other hand, $\frac{1}{n} h_K(u)f_K(u)d\s(u)$ is the volume element for $K$ and then \begin{equation*} h_{\phi K} (v)
f_{\phi K}(v)d\s(v)=|det (\phi)|h_K(u)f_K(u)d\s(u). 
\end{equation*} This further implies 
\begin{equation*} \prod_{i=1}^n [h_{\phi K_i} (v)
f_{\phi K_i}(v)]^{\frac{1}{n}}d\s(v)=|det (\phi)| \prod_{i=1}^n[h_{K_i}(u)f_{K_i}(u)]^{\frac{1}{n}} d\s(u).
\end{equation*}
Combining with formulas (\ref{p-volume:convex}) and (\ref{support:ratio}), one has for $L_1, \cdots, L_n \in
\cK_0$,  \begin{eqnarray*}
 V_{p}(\phi \bK; \phi \cL) &=&  \frac{1}{n} \int_{S^{n-1}}
  \prod _{i=1}^n \left[h
_{ \phi L_i} (v)^{p} f_p(\phi K_i, v)\right]^{\frac{1}{n}}\,d\s(v)\\
 &=& \frac{1}{n} \int_{S^{n-1}} \prod _{i=1}^n   \left(\frac{h
_{ \phi L_i} (v)}{h_{\phi K_i}(v)}\right)^{\frac{p}{n}}    \prod _{i=1}^n \left[h
_{ \phi K_i}(v)  f_{\phi K_i}(v)\right]^{\frac{1}{n}}\,d\s(v)\\&=& |det(\phi) |\cdot  \frac{1}{n} \int_{S^{n-1}} \prod _{i=1}^n   \left(\frac{h
_{ L_i} (u)}{h_{K_i}(u)}\right)^{\frac{p}{n}}    \prod _{i=1}^n \left[h
_{K_i}(u)  f_{K_i}(u)\right]^{\frac{1}{n}}\,d\s(u)\\
&=&|det(\phi) |V_{p}(\bK; \cL).
\end{eqnarray*}
Letting $L_1=\cdots =L_n=L$ in formula (\ref{affine:dual:mixed:volume}), one gets  for $p\geq 0$,
\begin{eqnarray*}
G_p^{(1)} (\phi K_1, \cdots,\phi  K_n)\!\!&=&\!\!\inf_{ L\in \cK_0 }
 \left\{n V_p(\phi \bK;  \phi L, \cdots, \phi L) ^{\frac{n}{n+p}}\
 |(\phi L) ^\circ|^{\frac{p}{n+p}}\right\} \\
 \!\!&=&\!\!|det(\phi)|^{\frac{n}{n+p}}|det(\phi)|^{-\frac{p}{n+p}}  \inf_{ L\in \cK_0  }
 \left\{n V_p( \bK;   L, \cdots,  L) ^{\frac{n}{n+p}}\
 |L ^\circ|^{\frac{p}{n+p}}\right\}\\
 \!\!&=&\!\! |det(\phi)|^{\frac{n-p}{n+p}}G_p^{(1)}( K_1, \cdots,  K_n).
\end{eqnarray*} Similarly, for $-n\neq p<0$, 
 \begin{eqnarray*}
G_p^{(1)} (\phi K_1, \cdots,\phi  K_n)\!\!&=&\!\!\sup_{ L\in \cK_0 }
 \left\{n V_p(\phi \bK;  \phi L, \cdots, \phi L) ^{\frac{n}{n+p}}\
 |(\phi L) ^\circ|^{\frac{p}{n+p}}\right\}  \\
 \!\!&=&\!\! |det(\phi)|^{\frac{n-p}{n+p}}G_p^{(1)}( K_1, \cdots,  K_n).
\end{eqnarray*} The proofs for $G_p^{(2)}(\phi K_1, \cdots,\phi  K_n)$ and $G_p^{(3)}(\phi K_1, \cdots,\phi  K_n)$ follow the same line.

 \bp \label{comparison of three} Let $K_1, \cdots, K_n\in \cF_0$. \\ (i). For $p\geq 0$, one has
 $$G_p^{(1)} (K_1, \cdots,  K_n)\geq G_p^{(2)} (K_1, \cdots,   K_n)\geq G_p^{(3)} (K_1, \cdots, K_n). $$
 (ii). For $-n< p<0$, one has $$G_p^{(1)} (K_1, \cdots,  K_n)\leq G_p^{(2)} (K_1, \cdots,   K_n)\leq G_p^{(3)} (K_1, \cdots, K_n). $$
 (iii). For $p<-n$, one has  $$G_p^{(1)} (K_1, \cdots,  K_n)\leq G_p^{(3)} (K_1, \cdots, K_n) \leq G_p^{(2)} (K_1, \cdots,   K_n). $$ \ep 
 \noindent \textbf{Proof.} 
 (i). It is clear for $p=0$. Let $p>0$, then $(L, \cdots, L)\in \cK_0^n$ and inequality (\ref{dual mixed volume}) imply  \begin{eqnarray*}  \inf_{ L  \in \cK_0}
 \left\{nV_p(\bK;  L, \cdots, L) ^{\frac{n}{n+p}}\
|L^\circ|^{\frac{p}{n+p}}\right\}  &\geq& \inf_{  \cL\in \cK_0^n }
\bigg\{n V_p(\bK;  \cL) ^{\frac{n}{n+p}}\ \prod_{i=1}^n
|L_i^\circ|^{\frac{p}{(n+p)n}}\bigg\} \\ &\geq& \inf_{ \cL\in \cK_0^n}
\left\{n V_p(\bK;  \cL) ^{\frac{n}{n+p}}\
\widetilde{V}(\cL^\circ )^{\frac{p}{n+p}}\right\}.
\end{eqnarray*} 
According to Definition \ref{equivalent:mixed
geominimal:surface:area-1}, one gets the desired conclusion.\vskip 2mm \noindent 
(ii). Let $p\in (-n, 0)$, then $\frac{p}{n+p}<0$. By $(L, \cdots, L)\in \cK_0^n$ and inequality (\ref{dual mixed volume}), one has
 \begin{eqnarray*}  \sup_{ L  \in \cK_0}
 \left\{nV_p(\bK;  L, \cdots, L) ^{\frac{n}{n+p}}\
|L^\circ|^{\frac{p}{n+p}}\right\}  &\leq& \sup_{  \cL\in \cK_0^n }
\bigg\{n V_p(\bK;  \cL) ^{\frac{n}{n+p}}\ \prod_{i=1}^n
|L_i^\circ|^{\frac{p}{(n+p)n}}\bigg\} \\ &\leq& \sup_{ \cL\in \cK_0^n}
\left\{n V_p(\bK;  \cL) ^{\frac{n}{n+p}}\
\widetilde{V}(\cL^\circ)^{\frac{p}{n+p}}\right\} .
\end{eqnarray*} According to Definition \ref{equivalent:mixed
geominimal:surface:area-1}, one gets the desired conclusion.\vskip 2mm \noindent
(iii). Let $p\in (-\infty, -n)$, then $\frac{p}{n+p}>0$. From inequality (\ref{dual mixed volume}), one has
\begin{eqnarray*}  \sup_{ \cL\in \cK_0^n}
\left\{n V_p(\bK;  \cL) ^{\frac{n}{n+p}}\
\widetilde{V}(\cL^\circ)^{\frac{p}{n+p}}\right\}  &\leq& \sup_{  \cL\in
\cK_0^n } \bigg\{n V_p(\bK;  \cL) ^{\frac{n}{n+p}}\
\prod_{i=1}^n |L_i^\circ|^{\frac{p}{(n+p)n}}\bigg\} ,
\end{eqnarray*}
which implies $G_p^{(3)} (K_1, \cdots,  K_n)\leq G_p^{(2)} (K_1,
\cdots,   K_n)$. Due to $(L,
\cdots, L)\in \cK_0^n$, one has
\begin{eqnarray*}  \sup_{ L  \in \cK_0}
 \left\{nV_p(\bK;  L, \cdots, L) ^{\frac{n}{n+p}}\
|L^\circ|^{\frac{p}{n+p}}\right\}  &\leq&  \sup_{ \cL\in \cK_0^n}
\left\{n V_p(\bK;  \cL) ^{\frac{n}{n+p}}\
\widetilde{V}(\cL^\circ)^{\frac{p}{n+p}}\right\} ,
\end{eqnarray*}
and hence  $G_p^{(1)} (K_1, \cdots,  K_n)\leq G_p^{(3)} (K_1,
\cdots,   K_n)$. The desired result then follows.

\bp\label{proposition 4.3}  Let $K_1, \cdots, K_n \in
\cF_0$. \\ (i). For $p\geq 0$, $$as_p (K_1,\cdots, K_n)\leq
G_p^{(\a)} (K_1,\cdots, K_n), \ \ \ \a=1, 2, 3.$$  (ii). For $-n< p<0$, $$as_p (K_1,\cdots,
K_n)\geq G_p^{(\a)}(K_1,\cdots, K_n), \ \ \ \a=1, 2, 3.$$ (iii). For $p<-n$,  one has 
\begin{eqnarray*} as_p(K_1,\cdots,
K_n)&\geq& G_p^{(\a)}(K_1,\cdots, K_n), \ \ \ \a=1,   3; \\
as_p^{(1)} (K_1,\cdots,
K_n)&\geq& G_p^{(\a)}(K_1,\cdots, K_n), \ \ \ \a=1, 2, 3.\end{eqnarray*}  \ep 
\vskip 2mm \noindent{\bf Proof.} (i). The case of $p=0$ holds trivially. We only prove the case $p>0$. By Proposition \ref{equal:formula},  Definition \ref{equivalent:mixed
affine:surface:area-1} and $\cK_0\subset\cS_0$, one can get 
\begin{eqnarray*}
as_p (K_1, \cdots, K_n)&=&\inf_{ \cL\in \cS_0^n }
\bigg\{n V_p(\bK;  \cL^\circ)
^{\frac{n}{n+p}}\ \widetilde{V}(\cL )^{\frac{p}{n+p}}\bigg\} \\
&\leq& \inf_{ \cL \in \cK_0^n } \bigg\{n V_p(\bK;  \cL^\circ)
^{\frac{n}{n+p}}\ \widetilde{V}(\cL )^{\frac{p}{n+p}}\bigg\}= G_p^{(3)}(K_1, \cdots, K_n).
\end{eqnarray*}
  Hence, Proposition \ref{comparison of three} implies that, for $K_1, \cdots, K_n\in \cF_0$ and for all $p>0$, 
 \begin{eqnarray*}
as_p (K_1, \cdots, K_n)\leq G_p^{(3)}(K_1, \cdots, K_n)\leq  G_p^{(2)}(K_1, \cdots, K_n)\leq G_p^{(1)}(K_1, \cdots, K_n).
\end{eqnarray*}
(ii). Let  $-n<p<0$.  Propositions  \ref{equal:formula} and  \ref{comparison of three},  Definition \ref{equivalent:mixed
affine:surface:area-1} and $\cK_0\subset\cS_0$ imply that    
\begin{eqnarray*}
as_p (K_1, \cdots, K_n)\!\!\!&=&\!\!\!\sup_{ \cL\in \cS_0^n }\!\!
\bigg\{n V_p(\bK;  \cL^\circ)
^{\frac{n}{n+p}}  \widetilde{V}(\cL )^{\frac{p}{n+p}}\!\bigg\} \geq  \sup_{ \cL \in \cK_0^n }\!\! \bigg\{n V_p(\bK;  \cL^\circ)
^{\frac{n}{n+p}}  \widetilde{V}(\cL )^{\frac{p}{n+p}}\!\bigg\}\\
\!\!\!&=&\!\!\! G_p^{(3)}(K_1, \cdots, K_n)
 \geq  G_p^{(2)}(K_1, \cdots, K_n)\geq G_p^{(1)}(K_1, \cdots, K_n).
\end{eqnarray*}
(iii).  Let $p<-n$. Propositions  \ref{equal:formula} and  \ref{comparison of three},  Definition \ref{equivalent:mixed
affine:surface:area-1} and $\cK_0\subset\cS_0$ imply that    
\begin{eqnarray*}
as_p (K_1, \cdots, K_n)\!\!\!&=&\!\!\!\sup_{ \cL\in \cS_0^n }\!\!
\bigg\{n V_p(\bK;  \cL^\circ)
^{\frac{n}{n+p}}  \widetilde{V}(\cL )^{\frac{p}{n+p}}\!\bigg\} \geq  \sup_{ \cL \in \cK_0^n }\!\! \bigg\{n V_p(\bK;  \cL^\circ)
^{\frac{n}{n+p}}  \widetilde{V}(\cL )^{\frac{p}{n+p}}\!\bigg\}\\
\!\!\!&=&\!\!\! G_p^{(3)}(K_1, \cdots, K_n)
 \geq G_p^{(1)}(K_1, \cdots, K_n).
\end{eqnarray*} Similarly, Proposition \ref{comparison of three},  Definition \ref{equivalent:mixed
affine:surface:area-1} and $\cK_0\subset\cS_0$ imply that \begin{eqnarray*} as_p^{(1)} (K_1, \cdots, K_n) &=&\sup_{ \cL \in \cS_0^n }
\bigg\{n  V_p(\bK;  \cL ^\circ)
^{\frac{n}{n+p}}  \prod_{i=1}^n |L_i|^{\frac{p}{(n+p)n}}\bigg\}\\ &\geq & \sup_{ \cL \in \cK_0^n }
\bigg\{n  V_p(\bK;  \cL ^\circ)
^{\frac{n}{n+p}}  \prod_{i=1}^n |L_i|^{\frac{p}{(n+p)n}}\bigg\}= G_p^{(2)}(K_1, \cdots, K_n)\\
  & \geq& \max\{G_p^{(1)}(K_1, \cdots, K_n), G_p^{(3)}(K_1, \cdots, K_n)\},
\end{eqnarray*} as desired. This concludes the proof of Proposition \ref{proposition 4.3}.

 Let $\cV_p$ for $-n\neq p\in \bbR$ be the subset of  $\cF_0^+$  \cite{Lu96, Y2} defined as
\begin{equation*}
\cV_p=\{K\in \cF_0^+: \exists Q \in \cK_0 \quad with \quad f_p(K,
u)=h_Q(u)^{-(n+p)}, \forall u\in S^{n-1}\}.
\end{equation*}
Note that $\cV_p\neq \emptyset$ for all $-n\neq p\in \bbR$ as $\ball \in \cV_p$. 
\bt Let $K _1, \cdots, K_n \in \cV_p$. Then for $-n\neq p\in \bbR$,
\begin{equation*}
G ^{(3)}_p (K_1, \cdots, K _n)= as_p (K_1, \cdots, K_n).
\end{equation*}
\et \vskip 2mm {\bf Remark.} It is easily checked that for all $\bK\in \cF_0^n$ and  for all $\alpha=1, 2, 3$,  $$G ^{(\alpha)} _0 (K_1, \cdots, K _n)= as _0 (K_1, \cdots,
K _n)=\int_{S^{n-1}} {\prod _{i=1}^n
[f_{K_i}(u)h_{K_i}(u)]^{\frac{1}{n}}} d\sigma(u).$$
 \noindent \textbf{Proof.}  Let $K_i \in \cV_p$ and  $p>0$. Proposition 3.3 in \cite{Y2} asserts that, for all $-n\neq p\in \bbR$, $ K_i \in
\cV_p$ implies $\L_p K_i\in \cK_0$. Thus, by Propositions \ref{affine:p:curvature} and \ref{proposition 4.3}, we have
\begin{eqnarray}
G_p^{(3)} (K_1, \cdots, K_n)&\geq& as _p (K_1, \cdots, K _n) \nonumber\\
&=& n \bigg(\frac{\o _n}{ |\L _p K_1|^\frac{1}{n} \cdots |\L _p
K_n|^\frac{1}{n}}\bigg) ^{\frac{n}{n+p}} \widetilde{V}(\L _p K_1,
\cdots, \L _p K_n) \nonumber \\
&=&  nV_p(\bK; (\L _p K_1) ^\circ, \cdots, (\L _p K_n)
^\circ) ^{\frac{n}{n+p}}
\widetilde{V}(\L _p K_1, \cdots, \L _p K_n)^{\frac{p}{n+p}} \nonumber \\
&\geq& \inf_{ \cL\in \cK_0^n } \bigg\{ nV_p(\bK; \cL^\circ) ^{\frac{n}{n+p}} \widetilde{V}(\cL )^{\frac{p}{n+p}}\bigg\}  = G _p ^{(3)} (K_1, \cdots, K _n). \nonumber
\end{eqnarray}
Hence, if $p>0$, $G ^{(3)}_p (K_1, \cdots, K _n)= as_p (K_1, \cdots, K_n)$ for 
all
$K_1, \cdots, K_n \in \cV_p$. Now let $ -n \neq p<0$. As above, one has \begin{eqnarray}
G _p ^{(3)} (K_1, \cdots, K _n) &\leq&  as _p (K_1, \cdots, K _n) \nonumber \\
&=&  nV_p(\bK; (\L _p K_1) ^\circ, \cdots, (\L _p K_n)
^\circ) ^{\frac{n}{n+p}}
\widetilde{V}(\L _p K_1, \cdots, \L _p K_n)^{\frac{p}{n+p}}  \nonumber \\
&\leq& \sup_{ \cL\in \cK_0^n } \bigg\{ nV_p(\bK; \cL^\circ) ^{\frac{n}{n+p}} \widetilde{V}(\cL)^{\frac{p}{n+p}}\bigg\} =G _p ^{(3)} (K_1, \cdots, K _n). \nonumber
\end{eqnarray}
Hence, if $-n\neq p<0$, $G ^{(3)}_p (K_1, \cdots, K _n)= as_p (K_1, \cdots, K_n)$ for
all $K_1, \cdots, K_n \in \cV_p$. 

\vskip 2mm Define the subset $\cV_{p,n}$ of  $(\cF_0^+)^n$ for $-n\neq p\in \bbR$  as
\begin{equation*}
\cV_{p, n}=\bigg\{\!\bK \in (\cF_0^+)^n: \exists Q \in \cK_0 \ s.t. \ \prod_{i=1}^n [f_p(K_i,
u)]^{\frac{1}{n}}=h_Q(u)^{-(n+p)}, \forall u\in S^{n-1}\!\bigg\}.
\end{equation*}

\bt Let $(K _1, \cdots, K_n) \in \cV_{p, n}$ with $p\neq 0, -n$. \\
(i). For $p>-n$ and $\alpha=1, 2, 3$,   $$G ^{(\alpha)}_p (K_1, \cdots, K _n)= as_p (K_1, \cdots, K_n).$$ (ii). For $p<-n$ and $\alpha=1, 3$,  $$G ^{(\alpha)}_p (K_1, \cdots, K _n)= as_p (K_1, \cdots, K_n).$$ 
\et
\noindent {\bf Proof.} Note that $(K _1,
\cdots, K_n) \in \cV_{p,n}$ for $-n\neq p\in \bbR$ implies $L_0\in \cK_0$ with
\begin{eqnarray*} \big[\r_ {L_0}(u)\big]^n= \bigg(\prod
_{i=1}^n f_p({K_i}, u)\bigg)^{\frac{1}{n+p}}, \ \  \forall u\in
S^{n-1}.
\end{eqnarray*}  \noindent (i). We first prove the case of $p>0$. Proposition  \ref{proposition 4.3} and formula (\ref{mixed p curvature-1-L}) imply
\begin{eqnarray*}
G_p^{(1)} (K_1, \cdots, K_n)&\geq& as _p (K_1, \cdots, K _n) = nV_p(\bK; L_0^\circ, \cdots, L_0^\circ) ^{\frac{n}{n+p}} |L_0|^{\frac{p}{n+p}}   \\
&\geq& \inf_{ L\in \cK_0 } \bigg\{ nV_p(\bK; L^\circ, \cdots, L^\circ) ^{\frac{n}{n+p}} |L|^{\frac{p}{n+p}} \bigg\}  = G _p ^{(1)} (K_1, \cdots, K _n),
\end{eqnarray*} where the second inequality follows from $L_0\in \cK_0$. 
Hence,  for $\bK \in \cV_{p, n}$ and for $p>0$,   $G ^{(1)}_p (K_1, \cdots, K _n)= as_p (K_1,
\cdots, K_n).$  Combining with Propositions \ref{comparison of three} and  \ref{proposition 4.3}, one has  \begin{eqnarray*} G_p^{(1)} (K_1, \cdots,  K_n)&\geq& G_p^{(2)} (K_1, \cdots,   K_n)\geq G_p^{(3)} (K_1, \cdots, K_n)\\ &\geq& as_p (K_1, \cdots,  K_n)= G_p^{(1)} (K_1, \cdots,   K_n), \end{eqnarray*}  as desired. Similarly for $ -n <p<0$, Propositions \ref{comparison of three} and  \ref{proposition 4.3}, and formula  (\ref{mixed p curvature-1-L}) imply  \begin{eqnarray*}
  G_p^{(1)} (K_1, \cdots,  K_n)&\leq& G_p^{(2)} (K_1, \cdots,   K_n)\leq G_p^{(3)} (K_1, \cdots, K_n)\\ &\leq& as _p (K_1, \cdots, K _n) = nV_p(\bK; L_0^\circ, \cdots, L_0^\circ) ^{\frac{n}{n+p}} |L_0|^{\frac{p}{n+p}}   \\
&\leq& \sup_{ L\in \cK_0 } \bigg\{ nV_p(\bK; L^\circ, \cdots, L^\circ) ^{\frac{n}{n+p}} |L|^{\frac{p}{n+p}} \bigg\}  = G _p ^{(1)} (K_1, \cdots, K _n). \end{eqnarray*} 
  \noindent (ii). Similarly for $p<-n$, Propositions \ref{comparison of three} and  \ref{proposition 4.3}, and formula  (\ref{mixed p curvature-1-L}) imply  \begin{eqnarray*}
  G_p^{(1)} (K_1, \cdots,  K_n)& \leq& G_p^{(3)} (K_1, \cdots, K_n)\leq  as _p (K_1, \cdots, K _n)\\ &= & nV_p(\bK; L_0^\circ, \cdots, L_0^\circ) ^{\frac{n}{n+p}} |L_0|^{\frac{p}{n+p}}   \\
&\leq& \sup_{ L\in \cK_0 } \bigg\{ nV_p(\bK; L^\circ, \cdots, L^\circ) ^{\frac{n}{n+p}} |L|^{\frac{p}{n+p}} \bigg\}  = G _p ^{(1)} (K_1, \cdots, K _n).
\end{eqnarray*}  Hence, for $p<-n$ and for $\alpha=1, 3$,  $G ^{(\alpha)}_p (K_1, \cdots, K _n)= as_p (K_1, \cdots, K_n)$ .

 \vskip 2mm \noindent {\bf Remark.} Note that if all $K_i=K\in \cF_0$,  Definitions \ref{equivalent:mixed
 geominimal:surface:area-1} and \ref{p-geominimal} together with formula (\ref{coincide:mixed}) imply    $G_p^{(1)} (K,\cdots,
 K) =\widetilde{G}_p(K)$ for all $ -n\neq p\in  \bbR$.    Moreover,  if  $K\in \cV_p$, then $(K, \cdots, K) \in \cV_{p, n}$.  Therefore, for $K\in \cV_p$, \begin{eqnarray*} G ^{(\alpha)}_p (K, \cdots, K)&=& \widetilde{G}_p (K)=as_p(K),  \ \ for\ -n\neq p\in \bbR \ \mbox{and $\alpha=1,  3$}; \\  G ^{(2)}_p (K, \cdots, K)&=& \widetilde{G}_p (K)=as_p(K), \ \ for \ p>-n.\end{eqnarray*}  In particular, as $(\ball, \cdots, \ball)\in \cV_{p, n}$, one gets, 
\begin{eqnarray*} G_p^{(\a)} (\ball,\cdots,
 \ball)&=&\widetilde{G}_p(\ball)=n |\ball|, \ \ \ for \ -n\neq p\in \bbR\ and \ \a=1, 3; \\G_p^{(2)}(\ball,\cdots,
  \ball)&=&\widetilde{G}_p(\ball)=n |\ball|, \ \ \ for \ p>-n.  \end{eqnarray*}

   \section{Inequalities for  $L_p$ mixed geominimal surface areas}\label{section:inequalities}
In this section, we will prove the Alexander-Fenchel type inequality, the affine isoperimetric inequality, the Santal\'{o} style inequality, and the cyclic inequality.
\subsection{The Alexander-Fenchel type
inequality} 
The classical Alexander-Fenchel inequality for the mixed volume  (see
 \cite{LSW, Schn}) is fundamental in applications. It has been extended to the mixed
$p$-affine surface area \cite{Lu75, Lut1987, Lu96, WY2}. Here, we prove the Alexander-Fenchel type
inequality for $L_p$ mixed geominimal surface areas.

\bt Let $K_{1}, \cdots, K_{n}\in \cF_{0} $. For $1 \leq m \leq
n$ and for $\a=1, 2 $, one has, 
\begin{eqnarray*}
\big[G_{p}^{(\a)}(K_{1},\cdots, K_{n})\big]^{m}&\leq& \prod_{i=0}^{m-1}G_{p}^{(\a)}(K_{1},
\cdots, K_{n-m}, \underbrace{K_{n-i},\cdots, K_{n-i}}_{m}), \ \ \ -n<p<0. 
\end{eqnarray*} In particular, if $m=n$, then  for $-n<p<0$ and for $\a=1, 2,$  \begin{eqnarray}
\big[G_{p}^{(\a)}(K_{1},\cdots, K_{n})\big]^{n} \leq  \widetilde{G}_{p} (K_{1})\cdots
\widetilde{G}_{p}  (K_{n}).\nonumber \end{eqnarray} 
  \et \vskip 0.2cm

 \noindent \textbf{Proof.}  By H\"{o}lder's inequality (see \cite{HLP}), one has 
\begin{eqnarray}
 V_p(\bK;  \cL)^m\!\!\!\!\!
 &=&\!\!\!\!\! \left(\frac{1}{n}\int_{S^{n-1}}H_{0}(u)H_{1}(u)\cdots H_{m}(u) d\s(u)\right)^m \nonumber \\
 \!\!\!\!\!&\leq&\!\!\!\!\! \prod_{i=0}^{m-1}\Big(\frac{1}{n}\int_{S^{n-1}}H_{0}(u)[H_{i+1}(u)]^md\s(u)\Big)  \nonumber \\
 \!\!\!\!\!&=&\!\!\!\!\!\prod_{i=0}^{m-1}\!V_{p} (K_{1},\!\cdots\!,\! K_{n-m}, \underbrace{K_{n-i},\!\cdots\!, K_{n-i}}_{m}; L_{1},\!\cdots\!, L_{n-m}, \underbrace{L_{n-i},\!\cdots\!, L_{n-i}}_{m}),\ \ \ \ \ \ \ \ \ \ \label{Holder:mixed:-001} 
\end{eqnarray} where for $i=0, \cdots, m-1,$ we let \begin{eqnarray*} H_{0}(u)&=&[h_{L_1}^{p}(u)f_{p}(K_{1}, u)
\cdots h_{L_{n-m}}^{p} (u)f_{p}(K_{n-m}, u)]^{\frac{1}{n}},\\ H_{i+1}(u)&=&[h_{L_{n-i}}^{p}(u)f_{p}(K_{n-i}, u)]^{\frac{1}{n}}. \end{eqnarray*}
Note that if $-n<p<0$, then $\frac{n}{n+p}>0$. Therefore,  
\begin{eqnarray}&& \big[G_{p}^{(2)}(K_{1},\cdots, K_{n})\big]^{m} =\sup_{\cL\in \cK_0^n} \big\{ nV_p(\bK;  \cL)^{\frac{n}{n+p}}\prod_{i=1}^n |L_i^\circ|^{\frac{p}{n(n+p)}}\big\}^m \nonumber\\
 &&\ \  \leq
 \prod_{i=0}^{m-1}\sup_{L_i \in \cK_0}  \big[nV_{p}(K_{1},\cdots, K_{n-m}, \underbrace{K_{n-i},\cdots, K_{n-i}}_{m}; L_{1},\cdots, L_{n-m}, \underbrace{L_{n-i},\cdots, L_{n-i}}_{m})^{\frac{n}{n+p}}\nonumber\\ & &\ \ \ \ \ \  \ \ \ \ \ \ \ \ \ \ \ \times |L_{n-i}^\circ|^{\frac{mp}{n(n+p)}} \prod_{i=1}^{n-m} |L_i^\circ|^{\frac{p}{n(n+p)}}\big]\nonumber\\ &&\ \  =\prod_{i=0}^{m-1}G_{p}^{(2)}(K_{1},
\cdots, K_{n-m}, \underbrace{K_{n-i},\cdots, K_{n-i}}_{m}).\label{A-F type -1}
\end{eqnarray}
The case of $\a=1$ follows directly by letting all $L_i=L$ in inequalities (\ref{Holder:mixed:-001}) and (\ref{A-F type -1}).

\bt Let $K_1, \cdots, K_n\in \cF_0$.  \\ (i). For $p\geq 0$,
\begin{eqnarray} \label{mixed geominimal-geominimal-inequality-2}
[G_{p}^{(3)}(K_{1},\cdots, K_{n})]^{n} &\leq& [G_{p}^{(2)}(K_{1},\cdots, K_{n})]^{n}\leq \widetilde{G}_{p} (K_{1})\cdots
\widetilde{G}_{p}  (K_{n}). \end{eqnarray}
(ii). For $p<-n$, \begin{eqnarray}
 [G_{p}^{(2)}(K_{1},\cdots, K_{n})]^{n} \geq  \widetilde{G}_{p} (K_{1})\cdots
\widetilde{G}_{p}  (K_{n}).  \label{mixed geominimal-geominimal-inequality-3}
\end{eqnarray}
\et {\bf Proof.} (i). Let $p\geq 0$, then $\frac{n}{n+p}>0$. Inequality (\ref{Holder:mixed:-001}) implies that \begin{eqnarray*}
 V_p(\bK;  \cL)^n&\leq &\prod_{i=1}^{n} V_{p} (K_{i},  L_{i}). 
\end{eqnarray*} Combining with Proposition \ref{comparison of three}, Definition \ref{equivalent:mixed
 geominimal:surface:area-1}, and Definition \ref{p-geominimal},  one has 
\begin{eqnarray}\big[G_{p}^{(3)}(K_{1},\cdots, K_{n})\big]^{n} &\leq & \big[G_{p}^{(2)}(K_{1},\cdots, K_{n})\big]^{n} \nonumber\\ &=&\inf_{\cL\in \cK_0^n} \big\{ nV_p(\bK;  \cL)^{\frac{n}{n+p}}\prod_{i=1}^n |L_i^\circ|^{\frac{p}{n(n+p)}}\big\}^n \nonumber\\
&\leq&
 \inf_{\cL \in \cK_0^n}   \prod_{i=1}^{n}\big[n V_{p} (K_{i},  L_{i})^{\frac{ n}{ n+p }}  |L_{i}^\circ|^{\frac{ p}{ n+p }}  \big]\nonumber\\ 
 &=&
 \prod_{i=1}^{n}\inf_{L_i \in \cK_0}  \big[n V_{p} (K_{i},  L_{i})^{\frac{ n}{ n+p }}  |L_{i}^\circ|^{\frac{ p}{ n+p }}  \big] \nonumber\\ &=&\widetilde{G}_{p} (K_{1})\cdots
\widetilde{G}_{p}  (K_{n}).\nonumber
\end{eqnarray}
(ii). Let $p<-n$, then  $\frac{n}{n+p}<0$. By Definition \ref{equivalent:mixed
 geominimal:surface:area-1}   and Definition \ref{p-geominimal},  one has 
\begin{eqnarray}   \big[G_{p}^{(2)}(K_{1},\cdots, K_{n})\big]^{n}\nonumber &=&\sup_{\cL\in \cK_0^n} \big\{ nV_p(\bK;  \cL)^{\frac{n}{n+p}}\prod_{i=1}^n |L_i^\circ|^{\frac{p}{n(n+p)}}\big\}^n \nonumber\\
 &\geq&
 \prod_{i=1}^{n}\sup_{L_i \in \cK_0}  \big[n V_{p} (K_{i},  L_{i})^{\frac{ n}{ n+p }}  |L_{i}^\circ|^{\frac{ p}{ n+p }}  \big]=\widetilde{G}_{p} (K_{1})\cdots
\widetilde{G}_{p}  (K_{n}).\nonumber
\end{eqnarray}
\subsection{The affine isoperimetric inequality } We will prove the affine isoperimetric inequality and the Santal\'{o} type inequality for the mixed $L_p$ geominimal surface areas. 

\bt\label{mixed-afiine-isoperimetric}
 Let $K_{1}, \cdots, K_n \in \cF_c$ be
convex bodies in $\cF_0$ with
 centroid at the origin.

\noindent (i). For $p\geq 0$ and $\a=2, 3$,
\begin{equation*}
\bigg( \frac{G_p ^{(\a)} (K_1,\cdots, K_n)}{G_p ^{(\a)}
(\ball,\cdots, \ball)}\bigg)^n \leq \min\bigg\{
\bigg(\frac{|K_{1}|}{|\ball|}\cdots
\frac{|K_{n}|}{|\ball|}\bigg)^{\frac{n-p}{n+p}}, \bigg(\frac{|K_{1}^\circ|}{|\ball|}\cdots
\frac{|K_{n}^\circ|}{|\ball|}\bigg)^{\frac{p-n}{n+p}}\bigg\},
\end{equation*} with equality if all $K_1,\cdots, K_n$ are ellipsoids that are dilates of one another.\\
 (ii). For $0\leq p\leq n$ and $\a=2, 3$,
\begin{equation*}
\bigg( \frac{G_p ^{(\a)} (K_1,\cdots, K_n)}{G_p ^{(\a)}
(\ball,\cdots, \ball)}\bigg)^{n+p}\leq \min\bigg\{  \frac{V(K_{1},\cdots,
K_{n})}{V(\ball,\cdots, \ball)},\  \frac{\widetilde{V}(\ball,\cdots, \ball)}{\widetilde{V}(K_{1}^\circ,\cdots,
K_{n}^\circ)}\bigg\}^{n-p},
\end{equation*} with equality if all $K_1,\cdots, K_n$ are ellipsoids that are dilates of one another.   
In particular,
$$
G_{n}^{(\a)}(K_{1},\cdots, K_{n})\leq G_{n} ^{(\a)} (\ball,\cdots,
\ball). $$   (iii). For $p> n$ and $\a=2, 3$,
$$
\bigg(\frac{G^{(\a)}_{p}(K_{1},\cdots,
K_{n})}{G^{(\a)}_{p}(\ball,\cdots, \ball)}\bigg)^{n+p}\leq
\min\bigg\{  \frac{V(K_{1}^\circ,\cdots,
K_{n}^\circ)}{V(\ball,\cdots, \ball)},\  \frac{\widetilde{V}(\ball,\cdots, \ball)}{\widetilde{V}(K_{1},\cdots,
K_{n})}\bigg\}^{p-n}, 
$$ with equality if all $K_1,\cdots, K_n$ are ellipsoids that are dilates of one another.\\
 (iv). For $p<-n$, $$ {G_{p} ^{(2)} (K_{1},\cdots, K_{n})} \geq n  \o_n^{\frac{2n}{n+p}}
 \big[{\widetilde{V}(K_{1}^\circ,\cdots,
K_{n}^\circ)} \big]^{\frac{p-n}{n+p}}.
$$ 
\et

 \noindent \textbf{Proof.} Note that, 
for all 
$p\geq 0$ and all $K\in \cF_c$, one has (see Theorem 4.2 in \cite{Y2}) \begin{equation}
\frac{\widetilde{G}_p(K)}{\widetilde{G}_p(\ball)} \leq\min\bigg\{\bigg(\frac{|K|}{|\ball|}\bigg)^{\frac{n-p}{n+p}},\ \  \bigg(\frac{|K^\circ|}{|\ball|}\bigg)^{\frac{p-n}{n+p}}\bigg\}.\label{Isoperimetric-p-Goeminimal}\end{equation} 
 (i). Recall that  $G_{p}^{(\a)}(\ball,\cdots,
 \ball)=\widetilde{G}_{p}(\ball)$ for all $\a=2,3$ and $p\geq 0$. By inequalities
 (\ref{mixed geominimal-geominimal-inequality-2})
 and  (\ref{Isoperimetric-p-Goeminimal}), one gets that for all $p\geq
 0$ and $\a=2, 3$, \begin{eqnarray}\bigg(\frac{G_{p}^{(\a)}(K_{1},\cdots,
K_{n})}{G_{p}^{(\a)}(\ball,\cdots, \ball)}\bigg)^n
\!\!\!\!&\leq&\!\!\!\! \frac{\widetilde{G}_{p}(K_{1})}{\widetilde{G}_{p}(\ball)}\cdots
\frac{\widetilde{G}_{p}(K_{n})}{\widetilde{G}_{p}(\ball)} \nonumber \\\!\!\!\! &\leq&\!\!\!\! \min\!\bigg\{\!
\bigg(\frac{|K_{1}|}{|\ball|}\cdots
\frac{|K_{n}|}{|\ball|}\bigg)^{\frac{n-p}{n+p}}\!, \bigg(\frac{|K_{1}^\circ|}{|\ball|}\cdots
\frac{|K_{n}^\circ|}{|\ball|}\bigg)^{\frac{p-n}{n+p}}\!\bigg\}.\label{Isoperimetric-mixed-Goeminimal}
 \end{eqnarray} Clearly, equality holds if all $K_1,\cdots, K_n$ are ellipsoids that are dilates of one another.  
 
 \noindent (ii). If
$0\leq p\leq n$, then $\frac{n-p}{n+p}\geq0$ and $\frac{p-n}{n+p}\leq 0$. Note that  $V(\ball,\cdots, \ball)=|\ball|$. Combining inequality (\ref{Minkowski inequality-mixed volume}) with inequality
(\ref{Isoperimetric-mixed-Goeminimal}), one has, for $\a=2, 3$ $$
\frac{G_{p} ^{(\a)} (K_{1},\cdots, K_{n})}{G_{p} ^{(\a)}
(\ball,\cdots, \ball)}\leq \left(\frac{V(K_{1},\cdots,
K_{n})}{V(\ball,\cdots, \ball)}\right)^{\frac{n-p}{n+p}}.
$$ Also note that $\widetilde{V}(\ball, \cdots, \ball)=|\ball|$. Combining  inequality (\ref{Isoperimetric-mixed-Goeminimal}) with inequality (\ref{dual mixed volume}), one gets 
for $\a=2, 3$ $$
\frac{G_{p} ^{(\a)} (K_{1},\cdots, K_{n})}{G_{p} ^{(\a)}
(\ball,\cdots, \ball)}\leq \bigg(\frac{\widetilde{V}(K_{1}^\circ,\cdots,
K_{n}^\circ)}{\widetilde{V}(\ball,\cdots, \ball)}\bigg)^{\frac{p-n}{n+p}}. 
$$ Clearly, equality holds if all $K_1,\cdots, K_n$ are ellipsoids that are dilates of one another. \\ (iii). If $p>n$, then $\frac{n-p}{n+p}<0$ and $\frac{p-n}{n+p}>0$. By inequalities (\ref{Isoperimetric-mixed-Goeminimal}) and (\ref{dual mixed volume}),   
for $\a=2, 3$  $$
\frac{G_{p} ^{(\a)} (K_{1},\cdots, K_{n})}{G_{p} ^{(\a)}
(\ball,\cdots, \ball)}\leq \bigg(\frac{\widetilde{V}(K_{1},\cdots,
K_{n})}{\widetilde{V}(\ball,\cdots,
\ball)}\bigg)^{\frac{n-p}{n+p}}.
$$ Combining inequality (\ref{Minkowski inequality-mixed volume}) with inequality
(\ref{Isoperimetric-mixed-Goeminimal}), one has, for $\a=2, 3$ $$
\frac{G_{p} ^{(\a)} (K_{1},\cdots, K_{n})}{G_{p} ^{(\a)}
(\ball,\cdots, \ball)}\leq \left(\frac{V(K_{1}^\circ,\cdots,
K_{n}^\circ)}{V(\ball,\cdots, \ball)}\right)^{\frac{p-n}{n+p}}.
$$ Clearly, equality holds if all $K_1,\cdots, K_n$ are ellipsoids that are dilates of one another.\\ 
 (iv). Note that $\widetilde{G}_p(\ball)=n|\ball|=n\o_n$. From  Theorem 4.2 in \cite{Y2}  and inequality
 (\ref{mixed geominimal-geominimal-inequality-3}), one gets that for all $p<-n$,
\begin{equation}
  \big[{G_{p}^{(2)}(K_{1},\cdots,
K_{n})} \big]^n 
\geq {\widetilde{G}_{p}(K_{1})} \cdots
 {\widetilde{G}_{p}(K_{n})} \geq
\bigg(\frac{|K_{1}^\circ|}{|\ball|}\cdots
\frac{|K_{n}^\circ|}{|\ball|}\bigg)^{\frac{p-n}{n+p}} (n\o_n)^n. \nonumber
\end{equation}
The desired result follows from inequality  (\ref{dual mixed
volume}). \vskip 2mm \noindent 
{\bf Remark.} When $\bK\in \cF_c^n$ with $V(K_1, \cdots, K_n)=|\ball|$ or with $\widetilde{V}(K_1^\circ, \cdots, K_n^\circ)=|\ball|$, then 
Theorem \ref{mixed-afiine-isoperimetric} implies that for all $p\in (0, n)$, $$G_p^{(\a)}(K_1, \cdots, K_n)\leq G_p^{(\a)}(\ball, \cdots,\ball), \ \ \a=2, 3.$$  That is, the mixed $L_p$ geominimal surface area $G_p^{(\a)}(K_1, \cdots, K_n)$ with $\a=2, 3$ attains the maximum at original-symmetric ellipsoids that are dilates of each other.  While for $p>n$, the mixed $L_p$ geominimal surface area $G_p^{(\a)}(K_1, \cdots, K_n)$ with $\a=2, 3$ attains its maximum at original-symmetric ellipsoids that are dilates of each other among $\bK\in \cF_c^n$ such that either $V(K_1^\circ, \cdots, K_n^\circ)=|\ball|$ or $\widetilde{V}(K_1, \cdots, K_n)=|\ball|$.  Although the condition $K_1, \cdots, K_n\in \cF_c$ is used in Theorems  \ref{mixed-afiine-isoperimetric} and 
 \ref{Theorem:Santalo}, and Corollary \ref{Corollary:5.1}, such a condition can be replaced by the following more general condition: some $K_i$ are in $\cF_c$ but others are in $\cF_s$.

\bc\label{Corollary:5.1}  Let $\mathcal {E}$ be an origin-symmetric ellipsoid and $K_1, \cdots, K_n\in \cF_c$.\\
 (i). For $0\leq p<n$ and $K_1, \cdots, K_n\subset\mathcal {E}$, one
has for $\a=2,3$
\begin{equation*}
G_{p}^{(\a)}(K_{1},\cdots, K_{n})\leq \widetilde{G}_{p}(\mathcal
{E}).
\end{equation*}
For $p=n$, the inequality holds for all $K_i\in
 \cF_c$ by (ii) of Theorem \ref{mixed-afiine-isoperimetric}.
\\
(ii). For $ p>n$ and $\mathcal {E}\subset K_1, \cdots, K_n$, one has for $\a=2,3$
\begin{equation*}
G_{p} ^{(\a)}(K_{1},\cdots, K_{n})\leq \widetilde{G}_{p}(\mathcal
{E}).
\end{equation*}
(iii). For $p<-n$ and $K_1, \cdots, K_n \subset\mathcal {E}$, one has
\begin{equation*}
G_{p}^{(2)}(K_{1},\cdots, K_{n})\geq \widetilde{G}_{p}(\mathcal
{E}).
\end{equation*}
 \ec
\noindent \textbf{Proof.} By Proposition \ref{proposition-1}, it is
enough to prove the proposition for $\mathcal {E}=\ball$. \\
(i). For $0\leq p < n$, one has $\frac{n-p}{n+p}>0$ and hence
$\big(\frac{|K_i|}{|\ball|}\big)^{\frac{n-p}{n+p}}\leq 1$ as
$K_i\subset\ball$. Combining with inequality (\ref{Isoperimetric-mixed-Goeminimal}), one has for
$\a=2,3$
\begin{equation*}
\bigg(\frac{G_{p}^{(\a)}(K_{1},\cdots,
K_{n})}{G_{p}^{(\a)}(\ball,\cdots, \ball)}\bigg)^n \leq
\bigg(\frac{|K_{1}|}{|\ball|}\cdots
\frac{|K_{n}|}{|\ball|}\bigg)^{\frac{n-p}{n+p}}\leq 1,
\end{equation*}
and the desired inequality holds.  \\
(ii). The condition $p>n$ implies $\frac{n-p}{n+p}<0$ and hence
$\left(\frac{|K_i|}{|\ball|}\right)^{\frac{n-p}{n+p}}\leq
1$ as $\ball\subset K_i$. Combining with inequality (\ref{Isoperimetric-mixed-Goeminimal}), one has for
$\a=2,3$
\begin{equation*}
\bigg(\frac{G_{p}^{(\a)}(K_{1},\cdots,
K_{n})}{G_{p}^{(\a)}(\ball,\cdots, \ball)}\bigg)^n \leq
\bigg(\frac{|K_{1}|}{|\ball|}\cdots
\frac{|K_{n}|}{|\ball|}\bigg)^{\frac{n-p}{n+p}}\leq 1,
\end{equation*}
and the desired inequality holds.   \\
 (iii).  For $p<-n$, by Corollary 4.1 in \cite{Y2} and  inequality (\ref{mixed geominimal-geominimal-inequality-3}),  one has, for $K_i \subset \ball$  

\begin{equation*}
G_{p}^{(2)}(K_{1},\cdots, K_{n}) \geq \big[ \widetilde{G}_{p} (K_{1})\cdots
\widetilde{G}_{p}  (K_{n})\big]^{1/n}\geq \widetilde{G}_{p} (\ball), 
\end{equation*}
and the desired inequality holds. 

The celebrated Blaschke-Santal\'{o} inequality states that, for all $K\in \cK_c$ (or $K\in \cK_s$),  $|K||K^\circ|\leq |\ball|^2$ with equality if and only if $K$ is an origin-symmetric ellipsoid. For the lower bound, Bourgain and Milman proved the following inverse
Santal\'{o} inequality \cite{BM} (see also \cite{GK2, Nazarov2012}):  there is a constant $c>0$, such that, $c^n|\ball|^2\leq |K||K^\circ|$
for all $K\in \cK_c$ (or $K\in \cK_s$). The next theorem provides a Santal\'{o} type
inequality for $L_p$ mixed Goeminimal surface areas. 
\bt \label{Theorem:Santalo} Let $K_{1}, \cdots, K_n \in \cF_c$. \\
(i). For $p\geq 0$ and $\a=2,3$,
\begin{equation*}
G_p ^{(\a)} (K_1,\cdots, K_n) G_p ^{(\a)} (K_1 ^\circ,\cdots, K_n
^\circ)\leq [G ^{(\a)}_p(\ball, \cdots, \ball)] ^2.
\end{equation*} Equality holds if $K_1, \cdots, K_n$ are ellipsoids that are dilates to each other. \\ 
(ii). For $p<-n$
\begin{equation*}
G_p ^{(2)} (K_1,\cdots, K_n) G_p ^{(2)} (K_1 ^\circ,\cdots, K_n
^\circ)\geq c^n [\widetilde{G }_p(\ball)] ^2.
\end{equation*}
where c is the universal constant from the Bourgain-Milman inverse
Santal\'{o} inequality. \et

\noindent \textbf{Proof.} (i). From inequality (\ref{mixed
geominimal-geominimal-inequality-2}), we have for all $p\geq0$ and
$\a=2,3$
\begin{equation*}
\big[G_p ^{(\a)} (K_1,\cdots, K_n)  G_p ^{(\a)} (K_1 ^\circ,\cdots, K_n
^\circ)\big]^n \leq \widetilde{G}_p(K_1)\widetilde{G}_p(K_1^\circ) \cdots
\widetilde{G}_p(K_n)\widetilde{G}_p(K_n^\circ).
\end{equation*}
   Theorem 4.1 in \cite{Y2} implies  that for  $p\geq 0$ and $\a=2,3$
\begin{equation*}
\big[G_p ^{(\a)} (K_1,\cdots, K_n)  G_p ^{(\a)} (K_1 ^\circ,\cdots, K_n
^\circ)\big]^n \leq \big[\widetilde{G}_p(\ball)\big]^{2n}=[G ^{(\a)}_p(\ball, \cdots, \ball)]
^{2n},
\end{equation*} as desired. Clearly, equality holds if all $K_1, \cdots, K_n$ are origin-symmetric ellipsoids that are dilates of each other. \\
(ii). From inequality (\ref{mixed geominimal-geominimal-inequality-3}), we have
for $p<-n$,
\begin{equation*}
\big[G_p ^{(2)} (K_1,\cdots, K_n)  G_p ^{(2)} (K_1 ^\circ,\cdots, K_n
^\circ)\big]^n \geq \widetilde{G}_p(K_1)\widetilde{G}_p(K_1^\circ) \cdots
\widetilde{G}_p(K_n)\widetilde{G}_p(K_n^\circ).
\end{equation*}
Theorem 4.1 in \cite{Y2} implies  that 
 for all $p<-n$
\begin{equation*}
G_p ^{(2)} (K_1,\cdots, K_n) G_p ^{(2)} (K_1 ^\circ,\cdots, K_n
^\circ)  \geq c^n \big[\widetilde{G}_p(\ball)\big]^{2}. 
\end{equation*}

Let $S_p(K)$ denote the $p$-surface area of $K$ for all $p\in \bbR$. That is, $S_p(K)=nV_p(K, \ball)$. In particular, $S_p(\ball)=n|\ball|={G}_p^{(\a)}(\ball, \cdots, \ball)$ for all $\a=1, 2, 3$ if $p>-n$, and for all $\a=1, 3$ if $p<-n$.  
\bc Let $K_1, \cdots, K_n\in \cF_0$.  \\ (i). For $p\geq 0$, one has $$\frac{G_p^{(\a)}(K_1, \cdots, K_n)}{G_p^{(\a)}(\ball, \cdots, \ball)}\leq \prod_{i=1}^n \bigg(\frac{S_p(K_i)}{S_p(\ball)}\bigg)^{\frac{1}{n+p}}, \ \ \a=1, 2, 3.$$ Equality holds if all $K_i$ are balls with center at the origin.\\
(ii). For $p\in (-\infty, -n)$, one has $$\frac{G_p^{(\a)}(K_1, \cdots, K_n)}{G_p^{(\a)}(\ball, \cdots, \ball)}\geq \prod_{i=1}^n \bigg(\frac{S_p(K_i)}{S_p(\ball)}\bigg)^{\frac{1}{n+p}}, \ \ \ \a=1,  3.$$ Equality holds if all $K_i$ are balls with center at the origin.  For $\a=2$, one has
$$\frac{G_p^{(2)}(K_1, \cdots, K_n)}{n|\ball|}\geq \prod_{i=1}^n \bigg(\frac{S_p(K_i)}{S_p(\ball)}\bigg)^{\frac{1}{n+p}}.$$

 \ec
{\bf Remark.} Combining with Proposition \ref{proposition 4.3}, one has 
\begin{eqnarray*} \frac{as_p (K_1, \cdots, K_n)}{as_p (\ball, \cdots, \ball)}&\leq& \prod_{i=1}^n \bigg(\frac{S_p(K_i)}{S_p(\ball)}\bigg)^{\frac{1}{n+p}}, \ \ for\ 
p\geq 0;\\
\frac{as_p^{(1)}(K_1, \cdots, K_n)}{as_p(\ball, \cdots, \ball)}&\geq& \frac{as_p  (K_1, \cdots, K_n)}{as_p(\ball, \cdots, \ball)}\geq \prod_{i=1}^n \bigg(\frac{S_p(K_i)}{S_p(\ball)}\bigg)^{\frac{1}{n+p}}, \ \ for \ p<-n.
\end{eqnarray*} 
{\bf Proof.} Let $p\geq 0$. Definition \ref{equivalent:mixed geominimal:surface:area-1}, inequality (\ref{Holder:mixed:-001}) and Proposition \ref{comparison of three} imply that \begin{eqnarray*}
G_p^{(3)}(K_1, \cdots, K_n)&\leq& G_p^{(2)}(K_1, \cdots, K_n)\leq 
G_p^{(1)}(K_1, \cdots, K_n)\\&=&\inf_{  L\in \cK_0  }
 \left\{n V_p(\bK;  L, \cdots, L) ^{\frac{n}{n+p}} 
 |L ^\circ|^{\frac{p}{n+p}}\right\}\\&\leq &  
  n V_p(\bK;  \ball, \cdots, \ball) ^{\frac{n}{n+p}} 
  |\ball|^{\frac{p}{n+p}} \\ &\leq &  n |\ball|^{\frac{p}{n+p}} \bigg[\prod_{i=1}^n   V_p(K_i, \ball)\bigg] ^{\frac{1}{n+p}}=(n |\ball|)^{\frac{p}{n+p}} \bigg[\prod_{i=1}^n   S_p(K_i)\bigg] ^{\frac{1}{n+p}}.
\end{eqnarray*}  The desired inequality follows from dividing by $G_p^{(\a)}(\ball, \cdots,\ball)=n|\ball|=S_p(\ball)$ (with $\a=1, 2, 3$). Equality holds if all $K_i$ are balls with center at the origin. Similarly, for $p<-n$,  \begin{eqnarray*}
G_p^{(3)}(K_1, \cdots, K_n) & \geq& G_p^{(1)}(K_1, \cdots, K_n)=\sup_{L\in \cK_0}
 \left\{n V_p(\bK;  L, \cdots, L) ^{\frac{n}{n+p}} 
 |L ^\circ|^{\frac{p}{n+p}}\right\}\\&\geq &  
  n V_p(\bK;  \ball, \cdots, \ball) ^{\frac{n}{n+p}} 
  |\ball|^{\frac{p}{n+p}}   \geq   n |\ball|^{\frac{p}{n+p}} \bigg[\prod_{i=1}^n   V_p(K_i, \ball)\bigg] ^{\frac{1}{n+p}}\\ &=&(n |\ball|)^{\frac{p}{n+p}} \bigg[\prod_{i=1}^n   S_p(K_i)\bigg] ^{\frac{1}{n+p}},
\end{eqnarray*} where the second inequality is due to $\frac{n}{n+p}<0$. Dividing by $G_p^{(\a)}(\ball, \cdots,\ball)$ (with $\a=1, 3$), we get the desired inequality. Equality holds if all $K_i$ are balls with center at the origin. For $\a=2$ and $p<-n$, 
the above inequality together with Proposition \ref{comparison of three} imply 
\begin{eqnarray*}
G_p^{(2)}(K_1, \cdots, K_n)\geq  ( n |\ball|)^{\frac{p}{n+p}} \bigg[\prod_{i=1}^n   S_p(K_i)\bigg] ^{\frac{1}{n+p}},
\end{eqnarray*} and the desired inequality follows by dividing by $n|\ball|=S_p(\ball)$.

\subsection{The cyclic inequality}
\noindent We now prove cyclic inequalities
for $L_p$ mixed geominimal surface areas.

\bt \label{cyclic inequality }
Let $K_{1}, \cdots, K_n\in \cF_{0}$ and $\a=1, 2, 3$. \\

\noindent (i). If $-n<t<0<r<s$ or $-n<s<0<r<t$, then
\begin{equation*} G_r ^{(\a)}(K_1,\cdots,
K_n) \leq [G_s ^{(\a)}(K_1,\cdots,
K_n)]^{\frac{(t-r)(n+s)}{(t-s)(n+r)}}[G_t ^{(\a)}(K_1,\cdots,
K_n)]^{\frac{(r-s)(n+t)}{(t-s)(n+r)}}.
\end{equation*} (ii). If $-n<t<r<s<0$ or $-n<s<r<t<0$, then
\begin{equation*} G_r ^{(\a)}(K_1,\cdots,
K_n) \leq [G_s ^{(\a)}(K_1,\cdots,
K_n)]^{\frac{(t-r)(n+s)}{(t-s)(n+r)}}[G_t ^{(\a)}(K_1,\cdots,
K_n)]^{\frac{(r-s)(n+t)}{(t-s)(n+r)}}.
\end{equation*}  
(iii). If $t<r<-n<s<0$ or $s<r<-n<t<0$, then
\begin{equation*} G_r ^{(\a)}(K_1,\cdots,
K_n) \geq [G_s ^{(\a)}(K_1,\cdots,
K_n)]^{\frac{(t-r)(n+s)}{(t-s)(n+r)}}[G_t ^{(\a)}(K_1,\cdots,
K_n)]^{\frac{(r-s)(n+t)}{(t-s)(n+r)}}.
\end{equation*}
\et

\noindent \textbf{Proof.} We only prove the case of $\a=2$ and the cases $\a=1, 3$ follow along the same lines. Let $\bK\in \cF_{0}^n$ and $\cL\in \cK_0^n$. Let $r, s, t\in \bbR$ such that $0< \frac{t-r}{t-s}<1$. 
By formula (\ref{p-volume}) and H\"{o}lder inequality
(see \cite{HLP}), 
 \begin{eqnarray} nV_r(\bK; \cL) \!\!\!\!&=& \!\!\!\!  \int_{S^{n-1}}  \prod_{i=1}^n [h^{r}
_{L_i}(u)f_r(K_i, u)]^{\frac{1}{n}}d\s(u)\nonumber\\
 \!\!\!\!&=& \!\!\!\!\int_{S^{n-1}}  \bigg(\prod_{i=1}^n [h^{s} _{L_i}(u)f_s(K_i,
u)]^{\frac{1}{n}}\bigg)^{\frac{t-r}{t-s}} \bigg(\prod_{i=1}^n
[h^{t}
_{L_i}(u)f_t(K_i, u)]^{\frac{1}{n}}\bigg)^{\frac{r-s}{t-s}}d\s(u)\nonumber\\
 \!\!\!\!&\leq& \!\!\!\!\bigg(\!\int_{S^{n-1}}  \prod_{i=1}^n [h^{s} _{L_i}(u)f_s(K_i,
u)]^{\frac{1}{n}}d\s(u)\!\!\bigg)^{\frac{t-r}{t-s}}\bigg(\!\int_{S^{n-1}}
\prod_{i=1}^n
[h^{t} _{L_i}(u)f_t(K_i, u)]^{\frac{1}{n}}d\s(u)\!\!\bigg)^{\frac{r-s}{t-s}}\nonumber\\
 \!\!\!\!&=& \!\!\!\!\left[nV_s(\bK; \cL)\right]^{\frac{t-r}{t-s}}\left[n V_t(\bK;
\cL)\right]^{\frac{r-s}{t-s}}\label{p-mixed volume-4}.
\end{eqnarray}  
(i). Suppose that $-n<t<0<r<s$, which implies $0<
\frac{t-r}{t-s}<1$. For simplicity, we let $\lambda=\frac{(r-s)(n+t)}{(t-s)(n+r)}$, and in this case $\lambda>0$.  We also let $\cL_1=(L_{11}, \cdots, L_{n1})\in \cK_0^n$. Then,  \begin{eqnarray*}  [G_t
^{(2)}(K_1,\cdots\!,\!
K_n)]^{\lambda} 
\!=\!\bigg\{\!\sup_{ \cL_{1}\in
\cK_0^n}\!\! \bigg[n V_t(\bK; \cL_{1}) ^{\frac{n}{n+t}}\!\prod_{i=1}^n |L_{i1}
 ^\circ|^{\frac{t}{n(n+t)}}\!\bigg]\!\bigg\}^{\lambda}\!\!\!=\!\!\!\sup_{\cL_{1}\in
\cK_0^n}\!\! \bigg [ n V_t(\bK;  \cL_{1})
^{\frac{n}{n+t}}\!\prod_{i=1}^n |L_{i1}
 ^\circ|^{\frac{t}{n(n+t)}}\!\bigg]^{\lambda} .
\end{eqnarray*}
 By inequality
(\ref{p-mixed volume-4}) and $\frac{n}{n+r}>0$, one has, for all
$\cL=(L_1, \cdots, L_n)\in\cK_0^n$,
\begin{eqnarray} G_r ^{(2)}(K_1,\cdots,
K_n)   
 \!\!\!\!&\leq&  \!\!\!\! n\big[V_r(\bK; \cL)\big]^{\frac{n}{n+r}} 
 \prod_{i=1}^n|L_i
 ^\circ|^{\frac{r}{n(n+r)}} \nonumber\\
  \!\!\!\!&\leq& \!\!\!\! \bigg[n V_s(\bK; \cL)^{\frac{n}{n+s}} 
 \prod_{i=1}^n|L_i
 ^\circ|^{\frac{s}{n(n+s)}}\bigg]^{1-\lambda} \bigg[n  V_t(\bK; \cL)^{\frac{n}{n+t}} 
 \prod_{i=1}^n|L_i
 ^\circ|^{\frac{t}{n(n+t)}}\bigg]^{\lambda}  \nonumber\\
 \!\!\!\!&\leq& \!\!\!\! \bigg[n V_s(\bK; \cL)^{\frac{n}{n+s}} 
 \prod_{i=1}^n|L_i
 ^\circ|^{\frac{s}{n(n+s)}}\bigg]^{1-\lambda} \!\!\! \sup_{ \cL_{1}\in \cK_0 ^n}\!\! \bigg[n  V_t(\bK;  \cL_{1}) ^{\frac{n}{n+t}} \prod_{i=1}^n|L_{i1}
 ^\circ|^{\frac{t}{n(n+t)}}\bigg]^{\lambda} \nonumber\\
 \!\!\!\!&=& \!\!\!\! \bigg[n V_s(\bK; \cL)^{\frac{n}{n+s}} 
 \prod_{i=1}^n|L_i
 ^\circ|^{\frac{s}{n(n+s)}}\bigg]^{1-\lambda} [G_t ^{(2)}(K_1,\cdots,
K_n)]^{\lambda}. \label{inf-1}
\end{eqnarray} Taking infimum over $\cL\in\cK_0^n$ in inequality (\ref{inf-1}), we
have
\begin{eqnarray*} G_r ^{(2)}(K_1,\cdots,
K_n) &\leq&\inf_{ \cL \in
\cK_0 ^n } \bigg[ n  V_s(\bK; \cL) ^{\frac{n}{n+s}} 
 \prod_{i=1}^n|L_i
 ^\circ|^{\frac{s}{n(n+s)}}\bigg]^{1-\lambda} [G_t ^{(2)}(K_1,\cdots,
K_n)]^{\lambda}\\
&=&[G_s ^{(2)}(K_1,\cdots, K_n)]^{1-\lambda}\ [G_t
^{(2)}(K_1,\cdots, K_n)]^{\lambda},
\end{eqnarray*} where the equality follows from $1-\lambda=\frac{(r-s)(n+t)}{(t-s)(n+r)}>0$ and then \begin{eqnarray*} [G_s ^{(2)}(K_1,\cdots\!,\!
K_n)]^{1\!-\!\lambda} 
\!\!=\!\bigg\{\!\!\inf_{\cL\in \cK_0^n}\!\!\bigg[nV_s(\bK; \cL) ^{\frac{n}{n+s}} \!\!
 \prod_{i=1}^n|L_i
 ^\circ|^{\frac{s}{n(n+s)}}\!\bigg]\!\bigg\}^{1\!-\!\lambda}\!\!\!\!\!=\!\!\inf_{\cL \in
\cK_0 ^n}\!\! \bigg[nV_s(\bK; \cL) ^{\frac{n}{n+s}}\!\!
 \prod_{i=1}^n|L_i
 ^\circ|^{\frac{s}{n(n+s)}}\!\bigg]^{1\!-\!\lambda}\!\!\!.
\end{eqnarray*}
The case of $-n\!<\!s<\!0<\!r<\!t$ follows immediately by switching the roles
of $t$ and $s$. 
\vskip.2cm\noindent (ii). Suppose that $-n<t<r<s<0$, which clearly
implies $0<\frac{t-r}{t-s}<1$. We let  
$\lambda=\frac{(r-s)(n+t)}{(t-s)(n+r)}$ and in this case $\lambda>0$.  Therefore, we have, 
\begin{eqnarray*}  [G_t ^{(2)}(K_1,\cdots,
K_n)]^{\lambda} 
\!=\!\bigg\{\! \sup_{\cL \in \cK_0 ^n}\!\!  \bigg[n  V_t(\bK; \cL)
^{\frac{n}{n+t}} 
 \prod_{i=1}^n|L_i
 ^\circ|^{\frac{t}{n(n+t)}}\bigg]\!\bigg\}^{\lambda}\!\!=\!\!\sup_{  \cL\in
\cK_0 ^n}\!\! \bigg[n V_t(\bK; \cL) ^{\frac{n}{n+t}}
 \prod_{i=1}^n|L_i
 ^\circ|^{\frac{t}{n(n+t)}}\bigg]^{\lambda}.\end{eqnarray*}
Similarly, due to  $1-\lambda=\frac{(t-r)(n+s)}{(t-s)(n+r)}>0$, one has 
\begin{eqnarray*}  [G_s ^{(2)}(K_1,\cdots,
K_n)]^{1-\lambda} 
 =\sup_{  \cL \in
\cK_0 ^n }\bigg [nV_s(\bK; \cL) ^{\frac{n}{n+s}}\
 \prod_{i=1}^n|L_i
 ^\circ|^{\frac{s}{n(n+s)}}\bigg]^{1-\lambda}.
\end{eqnarray*}
By inequality (\ref{p-mixed volume-4}) and $\frac{n}{n+r}>0$, we
have, for all $\cL\in\cK_0^n$, \begin{eqnarray*}  nV_r(\bK; \cL)^{\frac{n}{n+r}}\prod_{i=1}^n|L_i
^\circ|^{\frac{r}{n(n+r)}}\leq \bigg[nV_s(\bK; \cL)^{\frac{n}{n+s}}\prod_{i=1}^n|L_i
^\circ|^{\frac{s}{n(n+s)}}\bigg] ^{1-\lambda} \bigg[nV_t(\bK;
\cL)^{\frac{n}{n+t}}\prod_{i=1}^n|L_i
^\circ|^{\frac{t}{n(n+t)}}\bigg]^{\lambda}.
\end{eqnarray*}
The desired inequality follows by taking supremum over  $\cL\in\cK_0^n$.
The case of $-n<s<r<t<0$ follows immediately by switching the roles
of $t$ and $s$. 
\vskip.2cm \noindent (iii). Suppose that $t<r<-n<s<0$, which clearly
implies $0<\frac{t-r}{t-s}<1$.  Let 
$\lambda=\frac{(r-s)(n+t)}{(t-s)(n+r)}$ and in this case $\lambda>1$. Then, 
\begin{eqnarray*} \big[G_t ^{(2)}(K_1,\cdots,
K_n)\big]^{\lambda}\!=\!\bigg\{\!\sup_{\cL\in \cK_0^n
} \!\!\bigg[n V_t(\bK;  \cL) ^{\frac{n}{n+t}}\prod_{i=1}^n
 |L_i
 ^\circ|^{\frac{t}{n(n+t)}}\bigg]\!\bigg\}^{\lambda}\!\!=\!\!\sup_{ \cL\in
\cK_0 ^n }\!\! \bigg[n V_t(\bK;  \cL) ^{\frac{n}{n+t}}\prod_{i=1}^n
 |L_i
 ^\circ|^{\frac{t}{n(n+t)}}\!\bigg]^{\lambda}\!.
\end{eqnarray*}
Similarly, due to  $1-\lambda=\frac{(t-r)(n+s)}{(t-s)(n+r)}<0$, one has
\begin{eqnarray*} \big[G_s ^{(2)}\!(K_1,\!\cdots\!,\!
K_n)\big]^{1\!-\!\lambda}\!\!=\!\bigg\{\!\!\sup_{\cL_1\!\in\!
\cK_0^n}\!\! \bigg[nV_s(\bK;  \cL_1) ^{\frac{n}{n+s}}\!\!\prod_{i=1}^n |L_{i1}
 ^\circ|^{\frac{s}{n(n+s)}}\!\bigg]\!\bigg\}^{1\!-\!\lambda}\!\!\!\!=\!\!\inf_{\cL_1\!\in\!
\cK_0^n}\!\! \bigg[n V_s(\bK;  \cL_1)
^{\frac{n}{n+s}}\!\prod_{i=1}^n
 |L_{i1}
 ^\circ|^{\frac{s}{n(n+s)}}\!\bigg]^{1\!-\!\lambda}\!\!\!.
\end{eqnarray*}
 By inequality (\ref{p-mixed volume-4}) and $\frac{n}{n+r}<0$, we
have, for all $\cL \in\cK_0^n$,
\begin{eqnarray*} G_r ^{(2)}(K_1,\cdots,
K_n)  
&\geq& nV_r(\bK; \cL)^{\frac{n}{n+r}}\prod_{i=1}^n |L_i
^\circ|^{\frac{r}{n(n+r)}}\\
 &\geq& \bigg[n  V_s(\bK; \cL)^{\frac{n}{n+s}}\prod_{i=1}^n
 |L_i
 ^\circ|^{\frac{s}{n(n+s)}}\bigg]^{1-\lambda}  \bigg[n V_t(\bK; \cL)^{\frac{n}{n+t}}\prod_{i=1}^n
 |L_i
 ^\circ|^{\frac{t}{n(n+t)}}\bigg]^{\lambda}  \\
 &\geq&\bigg[n V_t(\bK; \cL)^{\frac{n}{n+t}}\prod_{i=1}^n
 |L_i
 ^\circ|^{\frac{t}{n(n+t)}}\bigg]^{\lambda} \inf_{ \cL_1\in \cK_0^n}\bigg[n  V_s(\bK;  \cL_1) ^{\frac{n}{n+s}}\prod_{i=1}^n
 |L_{i1}
 ^\circ|^{\frac{s}{n(n+s)}}\bigg]^{1-\lambda} 
 \\&=&
  \bigg[n V_t(\bK; \cL)^{\frac{n}{n+t}}\prod_{i=1}^n
 |L_i
 ^\circ|^{\frac{t}{n(n+t)}}\bigg]^{\lambda} \  [G_s ^{(2)}(K_1,\cdots,
 K_n)]^{1-\lambda}.
\end{eqnarray*} The desired inequality follows by taking supremum over  $\cL\in\cK_0^n$.
The case of $s<r<-n<0<t$ follows immediately by switching the roles
of $t$ and $s$. 
\vskip 2mm  \noindent{\bf Remark.} Note that the statement of Theorem
\ref{cyclic inequality } does not include the cases of $s = 0$ or $r
= 0$ or $t = 0$. However, from the proof of Theorem \ref{cyclic
inequality }, one can easily see that cyclic inequalities still hold
for (only) one of $r, s, t$ equal to $0$.

The monotonicity of
$\left(\frac{\widetilde{G}_p(K)}{n|K|}\right)^{\frac{n+p}{p}}$ was proved in \cite{Y2}. Here
we prove similar results for $L_p$ mixed geominimal surface areas.
\bt Let $K_1, \cdots, K_n\in \cF_0$ be such that $G_0 ^{(\a)}(K_1,\cdots,
K_n)>0$. The $L_p$ mixed geominimal surface areas  are monotone increasing in
the following sense: for  $0 < q < p$, or $-n < q < 0 < p$, or $-n <
q < p < 0$, or $q < p < -n$,
\begin{equation*} \bigg[\frac{G_q ^{(\a)}(K_1,\cdots,
K_n)}{G_0 ^{(\a)}(K_1,\cdots, K_n)}\bigg]^{\frac{n+q}{q}} \leq
\bigg[\frac{G_p ^{(\a)}(K_1,\cdots, K_n)}{G_0 ^{(\a)}(K_1,\cdots,
K_n)}\bigg]^{\frac{n+p}{p}},  \  \ for \ all \ \a=1, 2, 3.
\end{equation*}
\et

\noindent \textbf{Proof.} 
\noindent  {\it Case 1: $0<q<p$.} Employing (i) of
Theorem \ref{cyclic inequality } to $t=0, r=q$ and $s=p$, then
\begin{equation*}
G_q ^{(\a)}(K_1,\cdots, K_n) \leq \big[G_p ^{(\a)}(K_1,\cdots,
K_n)\big]^{\frac{q(n+p)}{p(n+q)}}\big[G_0 ^{(\a)}(K_1,\cdots,
K_n)\big]^{\frac{(p-q)n}{p(n+q)}}.
\end{equation*}
We divide both sides of the inequality by $G_0 ^{(\a)}(K_1,\cdots,
K_n)$ and get
\begin{equation*}
\frac{G_q ^{(\a)}(K_1,\cdots, K_n)}{G_0 ^{(\a)}(K_1,\cdots, K_n)}
\leq \bigg[\frac{G_p ^{(\a)}(K_1,\cdots, K_n)}{G_0
^{(\a)}(K_1,\cdots, K_n)}\bigg]^{\frac{q(n+p)}{p(n+q)}}.
\end{equation*} 
The desired inequality follows by taking the power $\frac{n+q}{q}>0$ from both sides. 

\vskip 2mm \noindent {\it Case 2: $-n<q<0<p$.} Employing (i) of
Theorem \ref{cyclic inequality } to $r=0, t=q$ and $s=p$, then
\begin{equation*}
G_0 ^{(\a)}(K_1,\cdots, K_n) \leq \big[G_p ^{(\a)}(K_1,\cdots,
K_n)\big]^{\frac{q(n+p)}{n(q-p)}}\big[G_q ^{(\a)}(K_1,\cdots,
K_n)\big]^{\frac{(n+q)p}{n(p-q)}}.
\end{equation*}
We divide both sides of the inequality by $G_0 ^{(\a)}(K_1,\cdots,
K_n)$ and get
\begin{equation*}
1 \leq \bigg[\frac{G_p ^{(\a)}(K_1,\cdots, K_n)}{G_0
^{(\a)}(K_1,\cdots,
K_n)}\bigg]^{\frac{q(n+p)}{n(q-p)}}\bigg[\frac{G_q
^{(\a)}(K_1,\cdots, K_n)}{G_0 ^{(\a)}(K_1,\cdots,
K_n)}\bigg]^{\frac{(n+q)p}{n(p-q)}}.
\end{equation*}
The desired inequality follows by taking the power $\frac{(q-p)n}{pq}>0$ from both sides.

 \vskip 2mm \noindent {\it Case 3: $-n<q<p<0$.} Employing (ii) of
Theorem \ref{cyclic inequality } to $s=0, t=q$ and $r=p$, then
\begin{equation*}
G_p ^{(\a)}(K_1,\cdots, K_n) \leq \big[G_0 ^{(\a)}(K_1,\cdots,
K_n)\big]^{\frac{n(q-p)}{q(n+p)}}\ \big[G_q ^{(\a)}(K_1,\cdots,
K_n)\big]^{\frac{p(n+q)}{q(n+p)}}.
\end{equation*}
Dividing both sides of the inequality by $[G_0 ^{(\a)}(K_1,\cdots,
K_n)]^{\frac{n(q-p)}{q(n+p)}}$ and as $-n<p<0$, we get, 
\begin{equation*}
\bigg[\frac{G_q ^{(\a)}(K_1,\cdots, K_n)}{G_0 ^{(\a)}(K_1,\cdots,
K_n)}\bigg]^{\frac{n+q}{q}} \leq \bigg[\frac{G_p
^{(\a)}(K_1,\cdots, K_n)}{G_0 ^{(\a)}(K_1,\cdots,
K_n)}\bigg]^{\frac{n+p}{p}}.
\end{equation*}    {\it Case 4: $q<p<-n$.} Employing (iii) of
Theorem \ref{cyclic inequality } to $s=0, t=q$ and $r=p$, then
\begin{equation*}
G_p ^{(\a)}(K_1,\cdots, K_n) \geq [G_0 ^{(\a)}(K_1,\cdots,
K_n)]^{\frac{n(q-p)}{q(n+p)}}[G_q ^{(\a)}(K_1,\cdots,
K_n)]^{\frac{p(n+q)}{q(n+p)}}.
\end{equation*}
Dividing both sides of the inequality by $[G_0 ^{(\a)}(K_1,\cdots,
K_n)]^{\frac{n(q-p)}{q(n+p)}}$ and as $p<-n$, we get,
\begin{equation*}
\bigg[\frac{G_q ^{(\a)}(K_1,\cdots, K_n)}{G_0 ^{(\a)}(K_1,\cdots,
K_n)}\bigg]^{\frac{n+q}{q}} \leq\bigg[\frac{G_p
^{(\a)}(K_1,\cdots, K_n)}{G_0 ^{(\a)}(K_1,\cdots,
K_n)}\bigg]^{\frac{n+p}{p}}.
\end{equation*}

 \section{The $i$-th mixed $L_p$ geominimal surface areas}\label{section:ith}
 This section dedicates to the $i$-th mixed $L_p$ geominimal surface areas, in particular, its related (affine) isoperimetric inequalities.  
 Let $K, L\in \cF_0^+$ and $Q_1, Q_2\in \cK_0$, we define $V_{p, i}(K, L; Q_1, Q_2)$ for all $i\in \bbR$ and all $p\in \bbR$ as 
 \be \label{ith-p-mixed-volume-1} n V_{p, i}(K, L;
Q_1, Q_2)=\int_{S^{n-1}}[h_{Q_1} ^p(u)f_{p}(K, u)]
^{\frac{n-i}{n}}[h_{Q_2}^p(u)f_{p}(L, u)] ^{\frac{i}{n}}d\s(u). \ee  When $Q_1, Q_2\in \cS_0$, we use the variation formula for $V_{p, i}(K, L; Q_1^\circ, Q_2^\circ)$ as 
 \be  n V_{p, i}(K, L; Q_1^\circ, Q_2^\circ)=\int_{S^{n-1}}[\rho_{Q_1} ^{-p}(u)f_{p}(K, u)]
^{\frac{n-i}{n}}[\rho_{Q_2}^{-p}(u)f_{p}(L, u)] ^{\frac{i}{n}}d\s(u). \nonumber \ee 
We also define $\widetilde{V}_i(Q_1, Q_2)$ for all $i\in \bbR$ as follows:  $$n\widetilde{V}_i(Q_1, Q_2)=\int_{S^{n-1}} [\r_{Q_1}(u)]^{n-i} [\r_{Q_2}(u)]^{i}\,d\s(u).$$ By H\"{o}lder's inequality (see
\cite{HLP}), one has,
\begin{eqnarray} \label{dual mixed volume-1}
[\widetilde{V}_i(Q_1, Q_2)]^n &\leq& |Q_1|^{n-i} |Q_2|^i, \ \  \ if
\ 0<i < n;\\
 \label{dual mixed volume-2} [\widetilde{V}_i(Q_1, Q_2)]^n &\geq&
|Q_1|^{n-i} |Q_2|^i,\ \  \ if \ i <0 \ or \ i> n.
\end{eqnarray}
Equality holds in each inequality if and only if $Q_1$ and $Q_2$ are
dilates of each other. Equality always holds in (\ref{dual mixed
volume-1}) and (\ref{dual mixed volume-2}) for $i=0$ or $i=n$. 

For $K, L\in \cF_0^+,$ the $i$-th mixed $p$-affine surface area \cite{Lut1987, WL1, WY2} can be formulated as \be as_{p, i}(K, L)=\int_{S^{n-1}}[f_p(K, u)]^{\frac{n-i}{n+p}}[f_p(L, u)]^{\frac{i}{n+p}}\,d\s(u), \ \  -n\neq p\in \bbR,\ \  i\in \bbR. \label{ith-p-mixed-volume-2} \ee The $i$-th mixed $p$-affine surface area contains many functionals as its special cases, such as the $L_p$ affine surface area and the $p$-surface area (i.e., $i=-p$ and $L=\ball$). Related properties and (affine) isoperimetric inequalities can be found in \cite{WY2}. The following proposition provides an equivalent formula for the $i$-th mixed $p$-affine surface area. 
\bp\label{mixed i th affine surface}  Let $K, L\in \cF_0^+$.\\ 
(i). For $p\geq 0$, 
 \begin{eqnarray*} as_{p,i} (K, L)&=&\inf_{ Q\in \cS_0 }
 \left\{n [V_{p,i}(K, L;  Q^\circ, Q^\circ)]^{\frac{n}{n+p}}
 |Q |^{\frac{p}{n+p}}\right\}  \\
 &=&\inf_{ \{Q_1, Q_2\in \cS_0\} }
 \bigg\{n[V_{p,i}(K, L; Q_1^\circ, Q_2^\circ)]^{\frac{n}{n+p}}
 \widetilde{V}_i(Q_1, Q_2)^{\frac{p}{n+p}}\bigg\}.
\end{eqnarray*}
(ii). For $-n\neq p<0$,  \begin{eqnarray*} as_{p,i} (K, L)&=&\sup_{ Q\in \cS_0 }
 \left\{n [V_{p,i}(K, L;  Q^\circ, Q^\circ)]^{\frac{n}{n+p}}
 |Q |^{\frac{p}{n+p}}\right\} \\ 
 &=&\sup_{ \{Q_1, Q_2\in \cS_0\} }
 \bigg\{n[V_{p,i}(K, L; Q_1^\circ, Q_2^\circ)]^{\frac{n}{n+p}}
 \widetilde{V}_i(Q_1, Q_2)^{\frac{p}{n+p}}\bigg\}.
\end{eqnarray*}
\ep
\noindent \textbf{Proof.} First, notice that for all $-n\neq p\in \bbR$ and $i\in \bbR$, we have 
\begin{eqnarray} as_{p, i}(K, L)= n [V_{p,i}(K, L;  Q_{0}^\circ,
Q_{0}^\circ)]^{\frac{n}{n+p}}
 [\widetilde{V}_i(Q_{0}, Q_{0})]^{\frac{p}{n+p}}, \label{definition:i-thmixedp----1}
\end{eqnarray}  where $Q_0\in \cS_0$ is defined by  
$\r_{Q_0}(u)=[f_p(K, u)^{\frac{n-i}{n}} f_p(L, u) 
^{\frac{i}{n}}]^{\frac{1}{n+p}}, \forall u\in S^{n-1}$.

\vskip 2mm \noindent (i). Clearly it holds for $p=0$. Let $p\in (0,\infty)$, then  $\frac{n}{n+p}\in (0,1)$. Employing H\"{o}lder inequality (see \cite{HLP}) to formula (\ref{ith-p-mixed-volume-2}),
 one has, for all $ K, L\in \cF_0^+$, for all $Q_1, Q_2 \in \cS_0$, and for all $i\in \bbR$,
\begin{eqnarray}
as_{p, i}(K, L) \!\!\!&=&\!\!\!\int_{S^{n-1}}\!\!\bigg(\!
[\r_{Q_1}(u)^{-p}f_p(K, u)]^{\frac{n-i}{n}} [\r_{Q_2}(u)^{-p}f_p(L,
u)]
^{\frac{i}{n}}\!\bigg)^{\frac{n}{n+p}}\bigg(\! \r_{Q_1}^{n-i}(u)\r_{Q_2}^i(u)\!\bigg)^{\frac{p}{n+p}} d\sigma(u)\nonumber  \\
\!\!\!&\leq&\!\!\!  n
[V_{p,i}(K, L;  Q_1^\circ, Q_2^\circ)]^{\frac{n}{n+p}}
 [\widetilde{V}_i(Q_1, Q_2)]^{\frac{p}{n+p}}.
\nonumber
\end{eqnarray} Taking infimum over $Q_1, Q_2\in\cS_0$ and  together with formula (\ref{definition:i-thmixedp----1}), one has
\begin{eqnarray*} as_{p, i}(K, L)&\leq& \inf_{ \{Q_1, Q_2\in\cS_0 \} } \left\{n [V_{p,i}(K, L;  Q_1^\circ, Q_2^\circ)]^{\frac{n}{n+p}}
 \widetilde{V}_{i}(Q_1, Q_2)^{\frac{p}{n+p}}\right\}\\ &\leq& \inf_{ Q\in \cS_0 }
 \bigg\{n[V_{p,i}(K, L; Q ^\circ, Q ^\circ)]^{\frac{n}{n+p}}
 |Q|^{\frac{p}{n+p}}\bigg\}\leq  as_{p, i}(K, L).
\end{eqnarray*}   

\noindent (ii). Note that $-n\neq p<0$
 implies $\frac{n}{n+p}>1$ or $\frac{n}{n+p}<0$. Employing H\"{o}lder inequality (see \cite{HLP}) to formula (\ref{ith-p-mixed-volume-2}),
 one has, for all $ K, L\in \cF_0^+$, for all $Q_1, Q_2 \in \cS_0$, and for all $i\in \bbR$,
\begin{eqnarray}
as_{p, i}(K, L) \!\!\!&=&\!\!\!\int_{S^{n-1}}\!\!\bigg(\!
[\r_{Q_1}(u)^{-p}f_p(K, u)]^{\frac{n-i}{n}} [\r_{Q_2}(u)^{-p}f_p(L,
u)]
^{\frac{i}{n}}\!\bigg)^{\frac{n}{n+p}}\bigg(\! \r_{Q_1}^{n-i}(u)\r_{Q_2}^i(u)\!\bigg)^{\frac{p}{n+p}} d\sigma(u)\nonumber  \\
\!\!\!&\geq&\!\!\!  n
[V_{p,i}(K, L;  Q_1^\circ, Q_2^\circ)]^{\frac{n}{n+p}}
 [\widetilde{V}_i(Q_1, Q_2)]^{\frac{p}{n+p}}.
\nonumber
\end{eqnarray} 
Taking supremum over $Q_1, Q_2\in\cS_0$ and  together with formula (\ref{definition:i-thmixedp----1}), one has 
\begin{eqnarray*} as_{p, i}(K, L)&\geq& \sup_{ \{Q_1, Q_2\in\cS_0 \} } \left\{n [V_{p,i}(K, L;  Q_1^\circ, Q_2^\circ)]^{\frac{n}{n+p}}
 \widetilde{V}_{i}(Q_1, Q_2)^{\frac{p}{n+p}}\right\}\\ &\geq& \sup_{ Q\in \cS_0 }
 \bigg\{n[V_{p,i}(K, L; Q ^\circ, Q ^\circ)]^{\frac{n}{n+p}}
 |Q|^{\frac{p}{n+p}}\bigg\}\geq  as_{p, i}(K, L).
\end{eqnarray*}    {\bf Remark.} By inequality (\ref{dual mixed volume-1}) and Proposition \ref{mixed i th affine surface}, one has, for
$p\geq 0$ and $0< i < n$,
 \begin{eqnarray*} as_{p,i} (K, L) &=& \inf_{ \{Q_1,
Q_2\in \cS_0 \}} \left\{n  V_{p,i}(K, L;  Q_1^\circ, Q_2^\circ)
^{\frac{n}{n+p}}\ \widetilde{V}_i(Q_1,Q_2)^{\frac{p}{n+p}}\right\}\\
&\leq& \inf_{ \{Q_1, Q_2\in \cS_0 \}} \left\{n  V_{p,i}(K, L;
Q_1^\circ, Q_2^\circ)
^{\frac{n}{n+p}}\ |Q_1|^{\frac{p(n-i)}{n(n+p)}}|Q_2|^{\frac{pi}{n(n+p)}}\right\}\\
&\leq& \inf_{ Q  \in \cS_0}
 \left\{nV_{p,i}(K, L; Q^\circ, Q ^\circ) ^{\frac{n}{n+p}}\
|Q|^{\frac{p}{n+p}}\right\} =as_{p,i} (K, L) .
\end{eqnarray*}  Similarly,  for $-n< p< 0$ with $0<i < n$, or  $p<-n$ with $i>n$
(or $i<0$), one has 
 \begin{eqnarray*} as_{p,i} (K, L)= \sup_{ \{Q_1, Q_2\in \cS_0\} }
 \bigg\{n [V_{p,i}(K, L; Q_1^\circ,  Q_2^\circ)] ^{\frac{n}{n+p}}
 |Q_1 |^{\frac{p(n-i)}{n(n+p)}}|Q_2  |^{\frac{pi}{n(n+p)}}\bigg\} .
\end{eqnarray*}
\vskip.2cm

Motivated by Proposition  \ref{mixed i th affine surface}, we define the $i$-th mixed $L_p$ geominimal surface areas as follows: \bd\label{ith-mixed
geominimal:surface:area-1}
Let $K,\ L\in \cF_0^+$, and  $\a=1,2,3$.   \\
(i). For $p=0$,  we let $$G_{0,i}^ {(\a)}(K, \ L)=\int_{S^{n-1}}
[h_K(u)f_K(u)]^{\frac{n-i}{n}} [h_L(u)f_
L(u)]^{\frac{i}{n}}d\sigma(u).$$ For $p>0$, 
 \begin{eqnarray*} G_{p,i}^{(1)}(K, L)&=&\inf_{ Q\in \cK_0 }
 \left\{n [V_{p,i}(K, L;  Q, Q)]^{\frac{n}{n+p}}
 |Q ^\circ|^{\frac{p}{n+p}}\right\};\\ 
 G_{p,i}^{(2)}(K, L)&=&\inf_{ \{Q_1, Q_2\in \cK_0\} }
 \bigg\{n [V_{p,i}(K, L; Q_1,  Q_2)] ^{\frac{n}{n+p}}
 |Q_1 ^\circ|^{\frac{p(n-i)}{n(n+p)}}|Q_2 ^\circ|^{\frac{pi}{n(n+p)}}\bigg\};\\
  G_{p,i}^{(3)}(K, L)&=&\inf_{ \{Q_1, Q_2\in \cK_0\} }
 \bigg\{n[V_{p,i}(K, L; Q_1, Q_2)]^{\frac{n}{n+p}}
 \widetilde{V}_i(Q_1^\circ, Q_2^\circ)^{\frac{p}{n+p}}\bigg\}.
\end{eqnarray*}
(ii). For $-n\neq p<0$,  
 \begin{eqnarray*} G_{p,i}^{(1)}(K, L)&=&\sup_{ Q\in \cK_0 }
 \left\{n [V_{p,i}(K, L;  Q, Q)]^{\frac{n}{n+p}}
 |Q ^\circ|^{\frac{p}{n+p}}\right\};\\ 
 G_{p,i}^{(2)}(K, L)&=&\sup_{ \{Q_1, Q_2\in \cK_0\} }
 \bigg\{n [V_{p,i}(K, L; Q_1,  Q_2)] ^{\frac{n}{n+p}}
 |Q_1 ^\circ|^{\frac{p(n-i)}{n(n+p)}}|Q_2 ^\circ|^{\frac{pi}{n(n+p)}}\bigg\};\\
  G_{p,i}^{(3)}(K, L)&=&\sup_{ \{Q_1, Q_2\in \cK_0\} }
 \bigg\{n[V_{p,i}(K, L; Q_1, Q_2)]^{\frac{n}{n+p}}
 \widetilde{V}_i(Q_1^\circ, Q_2^\circ)^{\frac{p}{n+p}}\bigg\}.
\end{eqnarray*}
\ed 
 Clearly, $G^{(\a)}_{p, i}(K, L)=G^{(\a)}_{p, n-i}(L, K)$ for all $i\in \bbR$. Moreover,
\be G^{(\a)}_{p, 0}(K, L)=\widetilde{G}_{p}(K), \ \ \& \ \
  G^{(\a)}_{p, n}(K, L)=\widetilde{G}_{p}(L).\label{equation-3} \ee  When $L=\ball$, we will write $G^{(\a)}_{p,
i}(K)$ as $G^{(\a)}_{p, i}(K, \ball)$ for all $\a=1,2,3$.

The $i$-th mixed $L_p$ geominimal surface areas have many properties similar to the mixed $L_p$ geominimal surface areas discussed in Section \ref{section:mixed}. For instance, the $i$-th mixed $L_p$ geominimal surface areas  are all affine invariant. Moreover, for $K, L\in \cF_0^+$, $i\in \bbR$, and $\a=1, 3$, one has \begin{eqnarray*}  G^{(\a)}_{p, i}(K, L) \geq as_{p, i}(K, L) \ \ if\  p\geq 0 \ \ and \ \  G^{(\a)}_{p, i}(K, L) \leq  as_{p, i}(K, L)\ \ if \  -n\neq p< 0.\end{eqnarray*} Equality holds if $K, L\in \cF_0^+$ satisfy the following property:  if $\exists Q \in \cK_0$  s.t. $$  [f_p(K ,
u)]^{\frac{n-i}{n}}[f_p(L,
u)]^{\frac{i}{n}}=h_Q(u)^{-(n+p)}, \ \ \forall u\in S^{n-1}. $$  In particular, one gets for all $i\in \bbR$, $-n\neq p\in \bbR$ and $\a=1, 3$,  \begin{eqnarray*}  G^{(\a)}_{p, i}(\ball, \ball )\!\!&=&\!\! as_{p, i}(\ball, \ball)=n|\ball|.\end{eqnarray*}
 
\bt \label{cyclic-theorem} Let $K, L\in\cF_0^+$. Suppose that $i, j, k\in
\mathbb{R}$ satisfy $ i< j<k$. Then, for $-n<p\leq 0$ and  $\a=1, 2$,  we have \be
\nonumber G ^{(\a)}_{p, j}(K, L)^{k-i}\leq G ^{(\a)}_{p, i}(K, L)^{k-j}G ^{(\a)}_{p,
k}(K, L)^{j-i}. \ee 
In particular, by letting $L=\ball$, we get
$$
  G ^{(\a)}_{p,
j}(K)^{k-i}\leq G ^{(\a)}_{p, i}(K)^{k-j}G ^{(\a)}_{p, k}(K)^{j-i}.
$$
  \et

 \noindent \textbf{Proof.} Let $K, L\in \cF_0^+$ and $Q_1, Q_2\in\cK_0$. From formula (\ref{ith-p-mixed-volume-1})  and
 H\"{o}lder's inequality, it follows that for all $p\neq-n$, and for $i< j<k$ (which clearly  implies
$0<\frac{k-j}{k-i}<1$),
 \begin{eqnarray}
V_{p, j}(K, L; Q_1, Q_2) \!\!\!&=&\!\!\! \frac{1}{n}\int_{S^{n-1}}[h_{Q_1}^ p(u)f_{p}(K,
u)] ^{\frac{n-j}{n}}[h_{Q_2}^ p (u)f_{p}(L, u)]
^{\frac{j}{n}}d\s(u)\nonumber\\
\!\!\!&\leq &\!\!\! \left\{\frac{1}{n}\int_{S^{n-1}}[h_{Q_1}^p(u)f_{p}(K,
u)]
^{\frac{n-i}{n}}[h_{Q_2}^p(u)f_{p}(L, u)] ^{\frac{i}{n}}d\s(u)\right\}^{\frac{k-j}{k-i}}\nonumber\\
&\quad&\times\left\{\frac{1}{n}\int_{S^{n-1}}[h_{Q_1 }^  p (u)f_{p}(K,
u)] ^{\frac{n-k}{n}}[h_{Q_2 }^ p(u)f_{p}(L, u)]
^{\frac{k}{n}}d\s(u)\right\}^{\frac{j-i}{k-i}}\nonumber\\
\!\!\! &=&\!\!\! V_{p, i}(K, L; Q_1, Q_2)^{\frac{k-j}{k-i}}V_{p, k}(K, L; Q_1,
Q_2)^{\frac{j-i}{k-i}}.  \label{p-mixed volume-5} 
\end{eqnarray}
 Note that $k-i>0, k-j>0$ and $j-i>0$. Then,  inequality (\ref{p-mixed volume-5}) implies that for $-n<p\leq 0$
\begin{eqnarray*}
G ^{(2)}_{p, j}(K, L)^{k-i}  &=&\sup_{ \{Q_1, Q_2\in \cK_0\} }
 \left\{n [ V_{p,j}(K, L;  Q_1, Q_2) ]^{\frac{n}{n+p}}\
 |Q_1 ^\circ|^{\frac{p(n-j)}{n(n+p)}}|Q_2 ^\circ|^{\frac{pj}{n(n+p)}}\right\}^{k-i}\\
 &\leq&\sup_{ \{Q_1, Q_2\in \cK_0\} }
 \left\{n [ V_{p,i}(K, L;  Q_1, Q_2) ]^{\frac{n}{n+p}}\
 |Q_1 ^\circ|^{\frac{p(n-i)}{n(n+p)}}|Q_2 ^\circ|^{\frac{pi}{n(n+p)}}\right\}^{k-j}\\
 &\quad&\times\sup_{ \{Q_1, Q_2\in \cK_0\} }
 \left\{n[V_{p,k}(K, L;  Q_1, Q_2) ]^{\frac{n}{n+p}}\
 |Q_1 ^\circ|^{\frac{p(n-k)}{n(n+p)}}|Q_2 ^\circ|^{\frac{pk}{n(n+p)}}\right\}^{j-i}\\
 &=&G ^{(2)}_{p, i}(K, L)^{k-j}G ^{(2)}_{p, k}(K, L)^{j-i}.
\end{eqnarray*}
The case $\a=1$ follows the same line by letting $Q_1=Q_2$.   
\vskip 2mm \noindent{\bf Remark.}  Let  $-n < p \leq 0$ and $\a=1,2$. For $0< i< n$, let $(i, j, k) = (0, i, n)$ in Theorem \ref{cyclic-theorem},  by formula (\ref{equation-3}), we have
 $$
  G ^{(\a)}_{p, i}(K, L)^{n}\leq G ^{(\a)}_{p, 0}(K, L)^{n-i}G ^{(\a)}_{p, n}(K, L)^{i}=\widetilde{G} _{p}(K)^{n-i}\widetilde{G} _{p}(L)^{i}.
 $$ The  equality always holds if  $i = 0$ or $i = n$.  Also note that for
 $-n<p\leq 0$, $\widetilde{G}_p(\ball)=n|\ball|$. By letting $L=\ball$, one has \begin{eqnarray*} 
 G ^{(\a)}_{p, i}(K)^{n}&\leq&
(n\omega_{n})^{i}\widetilde{G}_{p}(K)^{n-i}, \ \  0\leq i\leq 
n. \end{eqnarray*} Similarly,  for $i< 0$, let $(i, j, k) = (i, 0, n)$ in Theorem \ref{cyclic-theorem}, then
 $$
 G ^{(\a)}_{p, i}(K, L)^{n}G ^{(\a)}_{p, n}(K, L)^{-i}\geq G ^{(\a)}_{p, 0}(K, L)^{n-i} \Longleftrightarrow G ^{(\a)}_{p, i}(K, L)^{n}\geq \widetilde{G}_{p}(K)^{n-i}\widetilde{G}_{p}(L)^{i}.
 $$ For $i> n$, let $(i, j, k) = (0, n, i)$ in Theorem
\ref{cyclic-theorem}, then 
 $$
 G ^{(\a)}_{p, 0}(K, L)^{i-n}G ^{(\a)}_{p, i}(K, L)^{n}\geq G ^{(\a)}_{p, n}(K, L)^{i}\Longleftrightarrow G ^{(\a)}_{p, i}(K, L)^{n}\geq \widetilde{G}_{p}(K)^{n-i}\widetilde{G}_{p}(L)^{i}. $$   By letting $L=\ball$, we get \begin{eqnarray} 
G^{(\a)}_{p, i}(K)^{n}&\geq& (n\omega_{n})^{i}\widetilde{G}_{p}(K)^{n-i}, \ \  i<0 \ or\ i>n. \label{Minkowski-type-3}
 \end{eqnarray}

 \bt Let $K, L$ be convex bodies with continuous positive curvature functions and with centroid (or Santal\'{o} point) at the origin.  \\ (i). Let $p\geq 0$ and $0\leq i\leq n$, then for all $\a=2, 3$

$$
\frac{G^{(\a)}_{p,i}(K, L)}{n|\ball|}\leq
\min\bigg\{\! \!\!\left(\!\frac{|K|}{|\ball|}\!\right)^{\frac{(n-p)(n-i)}{n(n+p)}}\!\!\!\!\!\!, \ \ \left(\!\frac{|K^\circ|}{|\ball|}\!\right)^{\frac{(p-n)(n-i)}{n(n+p)}}\!\!\bigg\}\times \min\bigg\{\!\!\! \left(\!\frac{|L|}{|\ball|}\!\right)^{\frac{(n-p)i}{n(n+p)}}\!\!\!\!\!\!, \ \ \left(\frac{|L^\circ|}{|\ball|}\right)^{\frac{(p-n)i}{n(n+p)}}\bigg\}.$$   Moreover, 
$G ^{(\a)}_{p, i}(K, L)G ^{(\a)}_{p,
i}(K^{\circ}, L^\circ )\leq  (n\omega_{n})^{2}$. \\ (ii). Let $-n<p<0$ and $i\leq 0$, then for all $\a=1,2$
$$
\frac{G^{(\a)}_{p,i}(K)}{n|\ball|}\geq
\left(\frac{|K|}{|\ball|}\right)^{\frac{(n-p)(n-i)}{n(n+p)}}.
$$ Moreover, 
$G ^{(\a)}_{p, i}(K)G ^{(\a)}_{p,
i}(K^{\circ})\geq c^{n-i} (n\omega_{n})^{2}$ with $c>0$ a
universal constant. 
\\ (iii). Let $p<-n$ and $0\leq i\leq n$, then 
$$
\frac{G^{(2)}_{p,i}(K, L)}{n|\ball|}\geq
\left(\frac{|K^\circ|}{|\ball|}\right)^{\frac{(p-n)(n-i)}{n(n+p)}}\left(\frac{|L^\circ|}{|\ball|}\right)^{\frac{(p-n)i}{n(n+p)}}.
$$ Moreover, 
$G ^{(2)}_{p, i}(K, L)G ^{(2)}_{p,
i}(K^{\circ}, L^\circ)\geq c^{n} (n\omega_{n})^{2}$ with $c>0$ a universal constant.   \et
 \noindent \textbf{Proof.}  (i). The case of $p=0$ is clear. Let $p>0$ and $0\leq i \leq n$. H\"{o}lder's inequality implies 
\begin{eqnarray}
V_{p,i}(K,L; Q_1, Q_2)^n &=& \bigg(\int_{S^{n-1}}[h_{Q_1}^
p(u)f_{p}(K, u)] ^{\frac{n-i}{n}}[h_{Q_2}^ p (u)f_{p}(L, u)]
^{\frac{i}{n}}d\s(u)\bigg)^n   \nonumber \\
&\leq& \bigg(\int_{S^{n-1}}h_{Q_1}^ p(u)f_{p}(K, u)
d\s(u)\bigg)^{n-i}\bigg(\int_{S^{n-1}}h_{Q_2}^ p(u)f_{p}(L, u)
d\s(u)\bigg)^{i} \nonumber \\
 &=& V_p(K, Q_1)^{n-i}  V_p(L,
Q_2)^{i}. \label{p-mixed Minkowski}
\end{eqnarray}
Note that $\frac{n}{n+p}>0$.  Definitions \ref{ith-mixed geominimal:surface:area-1} and  
\ref{p-geominimal} together with inequality (\ref{p-mixed Minkowski}) imply \begin{eqnarray}
  [G^{(3)}_{p,i} (K, L)]^n   &\leq&  [G^{(2)}_{p,i} (K, L)]^n = \inf_{ \{Q_1, Q_2\in \cK_0\} }
 \bigg\{n [V_{p,i}(K, L; Q_1,  Q_2)] ^{\frac{n}{n+p}}
 |Q_1 ^\circ|^{\frac{p(n-i)}{n(n+p)}}|Q_2 ^\circ|^{\frac{pi}{n(n+p)}}\bigg\}^n\nonumber \\
&\leq& \inf_{ Q_1\in \cK_0 }\bigg\{n [V_p(K,
Q_1)]^{\frac{n}{n+p}}|Q_1 ^\circ|^{\frac{p}{n+p}}
\bigg\}^{n-i}\inf_{
 Q_2\in \cK_0 }\bigg\{n [V_p(K, Q_2)]^{\frac{n}{n+p}}|Q_2
^\circ|^{\frac{p}{n+p}}
\bigg\}^{i}\nonumber\\
&=& \widetilde{G}_p(K)^{n-i}  \widetilde{G}_p(L)^{i}.
\label{ith-minkowski-type-1}
 \end{eqnarray}
 Combining with inequalities (\ref{Isoperimetric-p-Goeminimal}) and 
(\ref{ith-minkowski-type-1}), one gets, for $0\leq i \leq n$, $p \geq 0$ and $\a=2, 3$, \begin{eqnarray*}
\frac{G^{(\a)}_{p,i}(K, L)}{n|\ball|}\!\!&\leq&\!\!
\left(\frac{\widetilde{G}_{p}(K)}{\widetilde{G}_{p}(\ball)}\right)^{\frac{n-i}{n}}\left(\frac{\widetilde{G}_{p}(L)}{\widetilde{G}_{p}(\ball)}\right)^{\frac{i}{n}}\\ \!\!&\leq&\!\!
\min\bigg\{\!\!\! \left(\!\frac{|K|}{|\ball|}\!\right)^{\frac{(n-p)(n-i)}{n(n+p)}}\!\!\!\!\!\!, \ \ \left(\!\frac{|K^\circ|}{|\ball|}\!\right)^{\frac{(p-n)(n-i)}{n(n+p)}}\!\!\bigg\}\times \min\bigg\{\!\!\! \left(\!\frac{|L|}{|\ball|}\!\right)^{\frac{(n-p)i}{n(n+p)}}\!\!\!\!\!\!, \ \ \left(\!\frac{|L^\circ|}{|\ball|}\!\right)^{\frac{(p-n)i}{n(n+p)}}\bigg\}.
\end{eqnarray*}
Recall that $
 \widetilde{ G}_{p}(K)\widetilde{G}_{p}(K^{\circ}) \leq (n\omega_{n})^{2}
 $ for all $K\in\cK_c$ and $p\geq0$  \cite{Y2}. Combining with inequality (\ref{ith-minkowski-type-1}), we have, for $0\leq i \leq n$, $p \geq 0$ and $\a=2, 3$,
 \begin{eqnarray*}
G ^{(\a)}_{p,i}(K, L)G ^{(\a)}_{p,i}(K^{\circ}, L^\circ) \leq
[\widetilde{G}_{p}(K)\widetilde{G}_{p}(K^{\circ})]^{\frac{n-i}{n}}[\widetilde{G}_{p}(L)\widetilde{G}_{p}(L^{\circ})]^{\frac{i}{n}}\leq   (n\omega_{n})^{2}.
\end{eqnarray*}

 \vskip 2mm \noindent (ii). Let $-n<p<0$ and $i\leq 0$.  Recall the affine isoperimetric inequality for the $L_p$ geominimal surface area in \cite{Y2} for $-n<p<0$ 
   \be
\nonumber \frac{\widetilde{G}_{p}(K)}{\widetilde{G}_{p}(\ball)}\geq
\left(\frac{|K|}{|\ball|}\right)^{\frac{n-p}{n+p}}. \ee
Combining with inequality (\ref{Minkowski-type-3}), one gets, for $i\leq 0$ and $\a=1, 2$, 
$$
\frac{G^{(\a)}_{p,i}(K)}{n|\ball|}\geq
\left(\frac{\widetilde{G}_{p}(K)}{\widetilde{G}_{p}(\ball)}\right)^{\frac{n-i}{n}}
\geq\left(\frac{|K|}{|\ball|}\right)^{\frac{(n-p)(n-i)}{n(n+p)}}.$$
Recall that $
 \widetilde{ G}_{p}(K)\widetilde{G}_{p}(K^{\circ}) \geq c^n (n\omega_{n})^{2} $ for all $K\in\cK_c$ and $-n<
p <0$ \cite{Y2}.  Combining with inequality (\ref{Minkowski-type-3}), we have, for $i\leq 0$ and $-n<p<0$ \begin{eqnarray*}
G ^{(\a)}_{p,i}(K)G ^{(\a)}_{p,i}(K^{\circ}) \geq
[\widetilde{G}_{p}(K)\widetilde{G}_{p}(K^{\circ})]^{\frac{n-i}{n}}[n\o_n ]^{\frac{2i}{n}}\geq  c^{(n-i)} (n\omega_{n})^{2}.
\end{eqnarray*}

 \noindent 
 (iii). Let $0\leq i\leq n$ and $p<-n$, which implies $\frac{n}{n+p}<0$. Combining with
 formula  (\ref{p-mixed Minkowski}), 
Definition \ref{ith-mixed geominimal:surface:area-1} and Definition
\ref{p-geominimal}, one has
\begin{eqnarray}
[G^{(2)}_{p,i} (K, L)]^n
 &=& \sup_{ \{Q_1, Q_2\in \cK_0\} }
 \bigg\{n [V_{p,i}(K, L; Q_1,  Q_2)] ^{\frac{n}{n+p}}
 |Q_1 ^\circ|^{\frac{p(n-i)}{n(n+p)}}|Q_2 ^\circ|^{\frac{pi}{n(n+p)}}\bigg\}^n\nonumber \\
&\geq& \sup_{ Q_1\in \cK_0 }\bigg\{n [V_p(K,
Q_1)]^{\frac{n}{n+p}}|Q_1 ^\circ|^{\frac{p}{n+p}}
\bigg\}^{n-i}\sup_{
 Q_2\in \cK_0 }\bigg\{n [V_p(K, Q_2)]^{\frac{n}{n+p}}|Q_2
^\circ|^{\frac{p}{n+p}}
\bigg\}^{i}\nonumber \\
&=& \widetilde{G}_p(K)^{n-i}  \widetilde{G}_p(L)^{i}. \label{ith-minkowski-type-2} \end{eqnarray}
Recall the affine isoperimetric inequality for the $L_p$ geominimal
surface area in \cite{Y2} for
 $p<-n$   \be
\nonumber \frac{\widetilde{G}_{p}(K)}{\widetilde{G}_{p}(\ball)}\geq
\left(\frac{|K^\circ|}{|\ball|}\right)^{\frac{p-n}{n+p}}. \ee
Combining with inequality  (\ref{ith-minkowski-type-2}), one has, for
$0\leq i\leq n$ and $p<-n$,
$$
\frac{G^{(2)}_{p,i}(K, L)}{n|\ball|}\geq
\left(\frac{\widetilde{G}_{p}(K)}{\widetilde{G}_{p}(\ball)}\right)^{\frac{n-i}{n}}\left(\frac{\widetilde{G}_{p}(K)}{\widetilde{G}_{p}(\ball)}\right)^{\frac{i}{n}}
\geq\left(\frac{|K^\circ|}{|\ball|}\right)^{\frac{(p-n)(n-i)}{n(n+p)}}\left(\frac{|L^\circ|}{|\ball|}\right)^{\frac{(p-n)i}{n(n+p)}}.
$$
 Recall that for all $K\in\cK_c$ and $p<-n$ \cite{Y2},   $
 \widetilde{ G}_{p}(K)\widetilde{G}_{p}(K^{\circ}) \geq c^n (n\omega_{n})^{2}
 $.  Combining with inequality (\ref{ith-minkowski-type-2}), we have, for $0\leq i \leq n$ and $p<-n$
 \begin{eqnarray*}
G ^{(2)}_{p,i}(K, L)G ^{(2)}_{p,i}(K^{\circ}, L^\circ) \geq
[\widetilde{G}_{p}(K)\widetilde{G}_{p}(K^{\circ})]^{\frac{n-i}{n}}[\widetilde{G}_{p}(L)\widetilde{G}_{p}(L^{\circ})]^{\frac{i}{n}}\geq  c^n (n\omega_{n})^{2}.
\end{eqnarray*}

\vskip 2mm \noindent {\bf Acknowledgments.} This research of DY is supported
 by a NSERC grant. The research of  BZ and JZ are supported by a NSFC grant (No.11271302).

\vskip 2mm \noindent Deping Ye, \ \ \ {\small \tt deping.ye@mun.ca}\\
{\small \em Department of Mathematics and Statistics\\
   Memorial University of Newfoundland\\
   St. John's, Newfoundland, Canada A1C 5S7 }

\vskip 2mm \noindent 
Baocheng Zhu, \ \ \ {\small \tt zhubaocheng814@163.com}\\
{\small \em School of Mathematics and Statistics\\
   Southwest University\\
    Chongqing 400715, China}\vskip 2mm \noindent Jiazu Zhou, \ \ \ {\small \tt zhoujz@swu.edu.cn}\\
{\small \em School of Mathematics and Statistics\\
   Southwest University\\
   Chongqing 400715, China }

\end{document}